\documentclass[11pt]{article}
\usepackage[francais]{babel}
\usepackage[latin1]{inputenc}

\usepackage{amsmath}
\usepackage{amsfonts}
\usepackage{amssymb}
\usepackage{amsthm}
\usepackage{graphicx,epsfig}
\usepackage{amscd}

\input{coupeq}

\newtheorem{Prop}{Proposition}[section]
\newtheorem{Def}[Prop]{Définition}
\newtheorem{Lem}[Prop]{Lemme}
\newtheorem{Thm}[Prop]{Théorème}
\newtheorem{Cor}[Prop]{Corollaire}
\newtheorem*{nThm}{Théorème}

\newcommand{\R}{\mathbb{R}}
\newcommand{\n}{\nabla}
\newcommand{\na}{\nabla ^{*}}
\newcommand{\tr}{\mathrm{tr\,}}
\newcommand{\im}{\mathrm{Im\,}}
\newcommand{\<}{\langle}
\renewcommand{\>}{\rangle}
\renewcommand{\d}{\partial}

\title{Rigidité infinitésimale de c\^ones-variétés Einstein à courbure
  négative}
\author{Grégoire Montcouquiol \\ \small{Laboratoire E.Picard (UMR 5580),
  Universit\'e Paul Sabatier, Toulouse}}

\begin{document}

\maketitle

\begin{minipage}{16cm}
\small

\noindent{\bf Abstract}
\vskip 0.5\baselineskip
\noindent

Starting with a compact hyperbolic cone-manifold of dimension $n \geq 3$, we
study the deformations of the metric with the aim of getting Einstein
cone-manifolds. If the singular locus is a closed codimension $2$ submanifold
and all cone angles are smaller than $2\pi$, we show that there is no
non-trivial infinitesimal Einstein deformations preserving the cone angles.

\vskip 0.5\baselineskip
\noindent{\bf R\'esum\'e}

\vskip 0.5\baselineskip
\noindent
Partant d'une cône-variété hyperbolique compacte de dimension $n \geq 3$, on
étudie les déformations de la métrique dans le but d'obtenir des
cônes-variétés Einstein. Dans le cas où le lieu singulier est une sous-variété
fermée de codimension $2$ et que tous les angles coniques sont plus petits que
$2\pi$, on montre qu'il n'existe pas de déformations Einstein infinitésimales
non triviales préservant les angles coniques.
\end{minipage}

\normalsize

\section{Introduction}

Dans leur célèbre article \cite{HK}, Hodgson et Kerckhoff montrent que pour
une large classe de cônes-variétés hyperboliques de dimension $3$, l'espace
des structures coniques hyperboliques au voisinage d'une cône-variété donnée
est paramétré par les angles coniques. Leur résultat principal est le théorème
de rigidité infinitésimale suivant : si $M$ est une cône-variété hyperbolique
de dimension $3$ de volume fini, dont le lieu singulier forme un entrelacs et
dont tous les angles coniques sont inférieurs à $2\pi$, alors il est
impossible de la déformer sans modifier ses angles.
Ce résultat, complété par des travaux plus récents
(cf notamment \cite{Kojima}, \cite{Weiss} et \cite{HK2}), a été le point de
départ de nombreux développements dans l'étude de la géométrie des variétés
hyperboliques de dimension $3$, tels que la géométrisation des petites
orbifolds ou l'étude des groupes kleiniens (\cite{BLP}, \cite{BB}).

Le principe de la démonstration du théorème de rigidité infinitésimale de
Hodgson et Kerckhoff est de
réussir à appliquer la méthode de Calabi-Weil aux
cônes-variétés : on montre que la représentation d'holonomie
n'admet pas de déformations non triviales de la forme voulue. Cela nécessite
d'établir des formules d'intégration par parties ainsi qu'un résultat du type
théorème de Hodge. Ce genre de difficultés est inhérent à l'étude
des cônes-variétés; nous les reverrons plus en détail.

Dans le cas des variétés fermées, Koiso \cite{Koiso} a donné un analogue de la
méthode de Calabi-Weil, qui n'utilise plus la représentation d'holonomie mais
étudie directement les déformations de la métrique (cf aussi \cite{Besse}, \S
12.H). Cette deuxième méthode
présente l'avantage d'être plus facilement généralisable et de s'appliquer, en
dimension supérieure, à une classe de variétés plus vaste, à savoir les
variétés Einstein (vérifiant de bonnes conditions de courbure).

Le but de ce papier est d'adapter la méthode de Koiso pour démontrer qu'en
dimension supérieure ou égale à trois, et sous des hypothèses voisines de
celles du théorème de Hodgson et Kerckhoff, on ne peut pas déformer une
cône-variété hyperbolique en des cônes-variétés Einstein sans en modifier les
angles coniques. En particulier, on redémontre dans le cas de la dimension
trois le théorème de rigidité infinitésimale ci-dessus.

\subsection{Présentation des résultats}

Le résultat que l'on se propose de démontrer ici est
le suivant :

\begin{nThm}[\ref{thmppal}]
Soit $M$ une c\^one-variété hyperbolique compacte, dont le lieu singulier
  forme une sous-variété fermée de codimension $2$, et dont tous les
  angles coniques sont strictement inférieurs à $2\pi$.
Alors toute déformation Einstein infinitésimale ne
  modifiant pas les angles coniques est triviale.
\end{nThm}

Les restrictions imposées à la géométrie des cônes-variétés sont
essentiellement les mêmes que dans l'article de Hodgson et
Kerckhoff \cite{HK}. On aurait pu remplacer l'hypothèse ``$M$ compacte''
par l'hypothèse ``$M$ de volume finie'', mais les choses sont quand même plus
simples dans le cas compact. La condition sur la géométrie du
lieu singulier est plus cruciale : c'est elle qui permet d'avoir
un bon modèle local et qui permet ainsi de faire les calculs. De
manière général, le lieu singulier d'une cône-variété peut être
beaucoup plus compliqué. Enfin, la condition sur les angles
coniques est une hypothèse technique qui paraît de prime abord
assez mystérieuse. En fait, on verra dans la section \ref{secLap}
que les angles coniques régissent en partie la croissance au voisinage du lieu
singulier des solutions d'un laplacien; plus les angles sont petits, plus on
contrôle ces solutions.

La définition précise des cônes-variétés envisagées se trouve dans
la section \ref{defnot}. Ce qu'il faut remarquer est que les
déformations infinitésimales d'une telle structure peuvent
toujours se mettre sous une forme standard (i.e. appartenant à une
certaine famille -- de dimension infinie -- de déformations) au
voisinage du lieu singulier. En particulier, une déformation ne
modifiant pas les angles a la propriété d'être $L^2$, à dérivée
covariante $L^2$. C'est entre autres pour cette raison que nous
travaillerons principalement dans le cadre $L^2$.

La section suivante rappelle la définition
des métriques et déformations infinitésimales Einstein; on y expose aussi le
problème des déformations triviales. Pour s'en débarrasser on cherche à
imposer la condition de jauge de Bianchi, ce qui revient à pouvoir résoudre
une équation de normalisation. On trouve ensuite dans la
section \ref{rapop} quelques résultats de la théorie des opérateurs non bornés
d'un espace de Hilbert qui nous serons utiles pous résoudre cette équation.

L'outil principal dans la démonstration de la rigidité infinitésimale est
connu sous le nom
de {\em technique de Bochner}. En partant d'une équation du type $Pu=0$ où $P$
est un opérateur différentiel du deuxième ordre de type Laplacien, on
exprime $P$ comme somme d'un opérateur auto-adjoint positif $Q^*Q$
de degré $2$ et d'un opérateur $R$ de degré $0$ faisant intervenir la
courbure. Une telle décomposition
$$P = Q^*Q + R$$
 s'appelle une {\em formule de Weitzenböck}; on en rencontrera à plusieurs
reprises dans ce texte. Ensuite, dans le cas où l'on travaille sur une
variété fermée, une intégration par parties donne
$$0 = \<Pu,u\> = ||Qu||^2 + \<Ru,u\>.$$
Si l'opérateur $R$ est tel que $\<Ru,u\> \geq c ||u||^2$ avec $c>0$, on
trouve alors $u=0$. Le lecteur intéressé par le sujet pourra se référer à
\cite{Besse}, \S 1.I.

Pour pouvoir adapter cette technique aux cônes-variétés il faut
pouvoir faire les intégrations par parties, et il est naturel pour
cela de travailler ici aussi avec des objets appartenant à des
espaces $L^2$. On donne dans la section \ref{secIPP} deux
résultats dans ce sens, ainsi que leur interprétation en termes
d'opérateurs non bornés.

La partie suivante (section \ref{secLap}) est le c\oe ur de ce texte. Elle
consiste en une étude détaillée de l'équation de normalisation et de
l'opérateur correspondant
$$L = \na \n + (n-1)Id = \Delta + 2(n-1)Id$$
agissant sur les 1-formes. Le but est de trouver des bons domaines sur
lesquels $L$ est auto-adjoint et donc inversible. Pour ce faire, et après
avoir préalablement exhibé une décomposition adaptée en séries de Fourier
généralisées (\S \ref{vois2}), on étudiera le comportement des
solutions de l'équation homogène au voisinage de la singularité. On
montrera que ce comportement est étroitement lié aux angles coniques; par
exemple, la norme ponctuelle d'une solution donnée au voisinage d'une
composante connexe du lieu singulier d'angle conique $\alpha$ est en $r^k$
avec $k \in \{\pm 1 \pm 2p\pi\alpha^{-1}, \pm 2p\pi\alpha^{-1} / p \in
\mathbb{Z} \}$. Les restrictions imposées sur les angles coniques permettent
de contrôler suffisamment les solutions de l'équation homogène, et finalement
les solutions de l'équation de normalisation tout court. On aboutit au
théorème suivant :

\begin{nThm}[\ref{nablad}] Soit $M$ une cône-variété hyperbolique dont tous les
  angles coniques sont strictement inférieurs à $2\pi$.
Soit $\phi$ une forme lisse appartenant à $L^2(T^*M)$. Alors il existe une
unique forme $\alpha \in
C^\infty(T^*M)$ solution de l'équation $Lu=\phi$ telle que $\alpha$, $\n
\alpha$, $d\delta \alpha$, et $\n d\alpha$ soient dans $L^2$.
\end{nThm}

Une fois ce résultat établi, il est relativement facile de faire fonctionner
la méthode de Koiso pour démontrer le théorème \ref{thmppal}; c'est
l'objet de la section \ref{rigidite}. Partant d'une
déformation infinitésimale Einstein $h_0$ préservant les angles
(donc à dérivée covariante $L^2$) d'une cône-variété hyperbolique,
dont tous les angles coniques sont inférieurs à $2\pi$,
la démonstration de sa trivialité se fait en deux temps. On a
d'abord besoin de se débarrasser des déformations triviales, on utilise donc
le résultat mentionné ci-dessus pour résoudre l'équation de normalisation.
On applique ensuite une technique de Bochner à la déformation normalisée
$h=h_0-\delta^*\alpha$. En utilisant la formule de Weitzenböck idoine et le
premier résultat d'intégration par parties, on obtient
$${\delta^\n d^\n h + (n-2) h =0}.$$
Une deuxième intégration par parties, un peu plus compliquée,
permet de conclure que $h_0 = \delta^*\alpha$, et donc que l'on a bien rigidité
infinitésimale relativement aux angles coniques au sein des cônes-variétés
Einstein.

\section{Les cônes-variétés et leurs déformations}\label{defnot}

Nous allons maintenant préciser le cadre dans lequel on se place.
La notion de c\^one-variété, plus générale que celle d'orbifold, a
été introduite par Thurston \cite{ThurstonGeom} pour l'étude des
déformations des variétés hyperboliques à cusps en dimension $3$.
Le cas le plus fréquemment rencontré est celui des
c\^ones-variétés à courbure constante. Celles-ci sont relativement
simples à définir, soit géométriquement comme un recollement de
simplexes géodésiques, soit en explicitant la métrique en
coordonnées; c'est cette dernière approche qui sera utilisée ici.
Le lecteur intéressé pourra se reporter à \cite{ThurstonShapes}
pour une définition par récurrence des cônes-variétés modelées sur
une géométrie.

La géométrie du lieu singulier d'une cône-variété arbitraire peut être très
compliquée. Dans le cadre de notre étude nous nous limiterons au cas où il
forme une sous-variété de codimension deux, ce qui permet de parler d'angle
conique le long de chaque composante connexe du lieu singulier et d'avoir des
bons modèles locaux pour mener à bien les calculs.

Enfin, comme notre but est de s'intéresser à des variétés Einstein, on
s'autorise une classe assez large de métriques à singularités : on demande
juste que la métrique conique ressemble asymptotiquement au produit de la
métrique du lieu singulier avec la métrique d'un cône (de dimension deux).

\bigskip

Soit $M$ une variété compacte de dimension $n\geq 3$, et $\Sigma =
\coprod_{i=1}^p \Sigma_i$ une sous-variété fermée plongée de
codimension $2$, dont les $\Sigma_i$ sont les composantes
connexes. Dans la suite de ce texte on emploiera souvent la
notation $M$ pour désigner improprement $M \setminus \Sigma$.

\begin{Def}\label{defcv} Soient $\alpha_1,\ldots,\alpha_p$ des réels
  positifs. La variété $M$ est munie d'une structure de {\em
  c\^{o}ne-variété}, de lieu singulier $\Sigma = \coprod_{i=1}^p \Sigma_i$ et
  d'angles coniques les $\alpha_i$, si :
\begin{itemize}

\item $M\setminus \Sigma$ est munie d'une métrique riemannienne $g$, non
  complète;

\item pour tout $i$, $\Sigma_i$ est munie d'une métrique riemannienne $g_i$;

\item pour tout $i$, tout point $x$ de $\Sigma_i$ a un voisinage $V$ dans
  $M$ difféomorphe à $D^2\times U$, avec $U=V\cap \Sigma_i$ un voisinage de
  $x$ dans $\Sigma_i$,
dans lequel $g$ s'exprime en coordonnées cylindriques locales sous
la forme $$g=dr^2 + (\frac{\alpha_i}{2\pi})^2 r^2d\theta^2+ g_i +
q,$$ où $q$ est un 2-tenseur symétrique vérifiant $g(q,q)=O(r^2)$
et $g(\n q,\n q)=O(r)$.
\end{itemize}
\end{Def}

\medskip

Dans la suite on exprimera souvent la métrique $g$ sous la forme
légèrement différente $$g = dr^2 + r^2d\theta^2+ g_i + q,$$
 où la coordonnée d'angle $\theta$ est définie non plus modulo
 $2\pi$ mais modulo l'angle conique $\alpha_i$.

Une c\^one-variété hyperbolique est alors une c\^one-variété telle
que les métriques $g$ et $g_i$ soient hyperboliques. On a dans ce
cas, en reprenant les notations de la définition, $$q =
(\frac{\alpha_i}{2\pi})^2(\sinh(r)^2 - r^2) d\theta^2 +
(\cosh(r)^2 - 1) g_i.$$

Pour démontrer la rigidité, nous aurons besoin que tous les angles coniques
soient inférieurs à $2\pi$, mais cette condition n'apparaît qu'à partir de la
fin de la partie \ref{secLap}.

\bigskip

Le caractère singulier des cônes-variétés pose certains problèmes pour adapter
la méthode de Koiso et faire fonctionner une technique de Bochner. Il faut
toujours vérifier si les choses marchent de la même manière que dans le cas
compact.

La première difficulté va venir des intégrations par parties. Premièrement,
pour garantir que les expressions manipulées ont un sens, nous serons obligés
de travailler avec des objets $L^2$. Deuxièmement, il va falloir démontrer
qu'on peut effectivement appliquer des formules de type Stokes : ce sera
l'objet de la partie \ref{secIPP}. Au final nous serons en mesure d'effectuer
des intégrations par parties pour les opérateurs $d$ et $\delta$, et $\n$ et
$\na$. Mais un tel résultat n'existe pas (à notre connaissance) pour les
opérateurs $d^\n$ et $\delta^\n$; nous devrons donc contourner cette
difficulté quand nous en aurons besoin (section \ref{rigidite}).

La plus grande difficulté va venir de l'équation de normalisation,
étudiée dans la section \ref{secLap}. Bien qu'en présence d'un
sympathique opérateur elliptique, on ne peut pas appliquer la
théorie classique sur une c\^one-variété, dont la métrique est
singulière. L'équation admettra encore des solutions, mais
celles-ci ne seront plus uniques, et on aura quoi qu'il arrive une
perte de régularité. Cependant, en imposant que les angles
coniques soient inférieurs à $2\pi$ nous aurons suffisamment de
contr\^ole sur la norme de certaines combinaisons linéaires des
dérivées des solutions pour faire fonctionner une technique de
Bochner.

\bigskip

Soit $(M,g)$ une cône-variété au sens ci-dessus, de lieu singulier
$\Sigma$. Soit maintenant $g_t$ une famille de métriques singulières,
dérivable, telle que $g_0=g$ et que pour tout $t$,
$(M,g_t)$ soit une cône-variété de lieu singulier $\Sigma$.

Si $x$ est un point de $\Sigma$, pour tout $t$ il existe par définition un
voisinage de $x$ dans $M$ dans lequel on a l'expression ci-dessus pour la
métrique en coordonnées cylindriques. Quitte à les restreindre, ces voisinages
sont tous difféomorphes, et on peut donc faire agir une famille $\phi_t$ de
difféomorphismes de telles façons que les coordonnées cylindriques locales
pour l'expression de $\phi_t^*g_t$ soient les mêmes pour tout $t$.

Dit d'une autre manière, il existe un voisinage $V$ de $x$ dans $M$,
difféomorphe à $D^2 \times U$ où $U = V \cap \Sigma$ est un voisinage de $x$
dans $\Sigma$, dans lequel on peut trouver des coordonnées cylindriques telles
que pour tout $t$, on ait  :
$$\phi_t^*g_t = dr^2 + (\frac{\alpha_t}{2\pi})^2r^2d\theta^2+ h_t + q_t.$$
Dans cette expression, $h_t$ désigne une métrique sur $U$ et $q_t$ est un
2-tenseur symétrique qui vérifie les conditions de la définition \ref{defcv}.

Finalement, quitte à modifier la famille $g_t$ par des difféomorphismes, ce
qui revient à modifier la déformation infinitésimale par une déformation
géométriquement triviale, on peut montrer que $h = \frac{dg_t}{dt}|_{t=0}$ est
au voisinage du lieu singulier combinaison linéaire des quatres types de
déformations suivants, modifiant :
\begin{list}{-}{}
\item l'angle,
\item la métrique du lieu singulier,
\item le reste,
\item et enfin, la façon de ``recoller'' la variable d'angle quand on passe
  d'un système de coordonnées à un autre.
\end{list}

Il est important de noter que les toutes ces déformations infinitésimales sont
$L^2$, mais que seules les trois dernières ont leur dérivée covariante dans
$L^2$. Ainsi, c'est au niveau du caractère $L^2$ ou non de la dérivée
covariante de la déformation que l'on voit si celle-ci préserve ou non les
angles coniques.

\section{Les métriques Einstein, leurs déformations et l'équation
de {normalisation}}\label{secEin}

Par définition, une {\em métrique Einstein} est une métrique
riemannienne $g$ vérifiant l'équation $$ric(g) = c g,$$ où le
terme de gauche est le tenseur de courbure de Ricci et où $c$ est
une constante. Notons que si on remplace $g$ par $\lambda g$, avec $\lambda$
une constante strictement positive, alors la nouvelle métrique vérifie
l'équation ci-dessus en remplaçant $c$ par $\lambda^{-1}c$;
donc en fait c'est principalement le signe et non la valeur exacte
de la constante $c$ qui compte. On peut ainsi distinguer trois
grandes classes de métriques Einstein suivant que $c$ est négatif,
positif ou nul.

Les métriques à courbure sectionnelle constante sont toujours
Einstein; en dimension $3$ ce sont les seules. Par contre dès la
dimension $4$ il y a beaucoup plus de métriques Einstein que de
métriques à courbure sectionnelle constante; on peut donc
considérer la condition Einstein comme un affaiblissement ou une
généralisation de la condition courbure sectionnelle constante.

Puisque l'on s'intéresse principalement aux cônes-variétés
hyperboliques, on ne considèrera que des métriques Einstein
vérifiant $E(g)=0$, avec $$E(g) = ric(g) + (n-1) g.$$ La constante
$(n-1)$ est choisie de telle sorte que les métriques hyperboliques
vérifient cette équation.

\bigskip

Soit $g_t$ une famille lisse de métriques Einstein (c'est-à-dire
vérifiant $E(g_t)= 0$) sur une variété donnée $M$, avec $g_0=g$.
Le $2$-tenseur symétrique $h=\frac{d}{dt}g_t|_{t=0}$ vérifie alors
l'équation d'Einstein linéarisée $$E'_g(h)=0.$$ Le calcul de
$E'_g$ est classique, cf par exemple \cite{Besse} \S 1.K :
$$E'_g(h) = \na_g\n_g h - 2 \mathring{R}_gh - \delta^*_g(2\delta_g
h + d\tr_g h).$$

Les opérateurs utilisés ici nécessitent un peu d'explication. La
notation $\n_g$, ou $\n$ pour simplifier, désigne la dérivée covariante ou
connexion de
Levi-Cività associée à la métrique riemannienne $g$. Elle admet un
adjoint formel noté $\na_g$ : si $(e_i)_{i=1\ldots n}$ est une base
orthonormée, on a $$\na_g \alpha (X_1,\ldots,X_p) = - \sum_{i=1}^n
(\n_{e_i} \alpha)(e_i,X_1,\ldots,X_p).$$

Pour les tenseurs symétriques, on définit $\delta^*_g :
\mathcal{S}^p M \to \mathcal{S}^{p+1} M$ comme étant la composée
de la dérivée covariante et de la symétrisation. En particulier,
si $\alpha \in \Omega^1 M = \mathcal{S}^1 M$, alors
\begin{eqnarray*}
\delta^*_g \alpha (x,y) & = &\frac{1}{2}((\n_x \alpha)(y) + (\n_y
\alpha)(x)\\
& = &\frac{1}{2}(g(\n_x
\alpha^\sharp,y) + g(\n_y \alpha^\sharp,x))\\
& = &\frac{1}{2}L_{\alpha^\sharp}g(x,y),
\end{eqnarray*}
où $L_{\alpha^\sharp}$ désigne la dérivée de Lie le long du champ
de vecteur $\alpha^\sharp$ dual (pour la métrique $g$) à la forme
$\alpha$. L'adjoint
formel de l'opérateur $\delta^*_g$ se note $\delta_g$ ; c'est juste la
restriction de $\na_g$ à $S^{p+1}M$.

Ensuite, $\mathring{R}_g$ désigne l'action du tenseur de courbure
$R_g$ sur les $2$-tenseurs symétriques : si $h$ est une section de
$S^2M$, on pose
$$\mathring{R}_g h(x,y) = \sum_{i=1}^n h(R_g(x,e_i)y,e_i),$$
où $(e_i)$ est une base orthonormale pour $TM$;
c'est encore un $2$-tenseur symétrique. Si $g$ est hyperbolique,
on a alors $$\mathring{R}_g h = h - (\tr_g h)g.$$

Enfin, la notation $\tr_g$ désigne juste la trace par rapport à
$g$ : si $h$ est un $2$-tenseur,
$$\tr_g h = \sum_{i=1}^n h(e_i,e_i).$$
Dans la suite et pour alléger les notations, on
omettra le plus fréquemment l'indice $g$.

\bigskip

Par définition, une {\em déformation Einstein infinitésimale} de
la variété Einstein $(M,g)$ est un $2$-tenseur symétrique $h$
vérifiant l'équation $E'_g(h)=0$.

Maintenant, si $g$ est Einstein et si $\phi$ est un
difféomorphisme de $M$, alors la métrique tirée en arrière $\phi^*
g$ est aussi Einstein. Par conséquent, si $\phi_t$ est une famille
lisse de difféomorphismes telle que $\phi_0$ soit l'identité,
alors la déformation infinitésimale associée
$\frac{d}{dt}\phi_t^*g|_{t=0}$ est naturellement Einstein. Une
telle déformation est qualifiée de {\em triviale}. Soit $X$ le
champ de vecteurs sur $M$ défini par
$X(x)=\frac{d}{dt}(\phi_t(x))|_{t=0}$, et soit $\alpha = X^\flat$
la $1$-forme duale, c'est-à-dire vérifiant $\alpha(Y) = g(X,Y)$
pour tout vecteur $Y$. On a les relations $$\frac{d}{dt}(\phi_t^*
g)|_{t=0} = L_X g = 2 \delta^*_g \alpha;$$ l'espace des
déformations infinitésimales triviales est donc égal à $\im
\delta^*_g$.

\bigskip

La façon habituelle de se débarrasser des déformations triviales
est d'imposer une condition de jauge, c'est-à-dire de ne
considérer que des déformations infinitésimales vérifiant une
certaine équation. On en trouve plusieurs dans la littérature, on
utilisera ici la jauge de Bianchi (cf \cite{Biquard} \S I.1.C,
\cite{Anderson} \S 2.3, à comparer à \cite{Besse} \S 12.C). On
veut donc que nos déformations infinitésimales $h$ vérifient
$$\beta_g(h)=0,$$ où $\beta_g : \mathcal{S}^2M \to \Omega^1M$ est
l'opérateur de Bianchi (associé à la métrique $g$) défini par
$$\beta_g(h)=\delta_g h +\frac{1}{2}d\tr_g h.$$

Ainsi, étant donnée une déformation infinitésimale $h_0$, on veut
pouvoir la modifier par une dé\-for\-ma\-tion triviale, de façon
essentiellement unique, de telle sorte que le résultat vérifie la
condition de jauge. Dit plus précisément, on veut trouver une
$1$-forme $\alpha$ telle que la déformation normalisée $h=h_0 -
\delta^*\alpha$ satisfasse $\beta(h)=0$; de façon équivalente, on
cherche à résoudre {\em l'équation de normalisation} (on omet les
indices) $$\beta\circ \delta^* \alpha = \beta(h_0).$$

L'étude de cette équation et de l'opérateur $\beta\circ\delta^*$ est l'objet
de la section \ref{secLap}. On se placera entre autre dans le cadre de la
théorie des opérateurs non bornés entre espace de Hilbert, dont les résultats
principaux sont cités dans la section suivante.

\section{Quelques rappels sur les opérateurs non
bornés}\label{rapop}

Nous allons annoncer un certain nombre de définitions et propriétés concernant
les opérateurs non bornés; le lecteur intéressé pourra consulter \cite{Riesz},
chapitre 8, ou \cite{Rudin}, chapitre 13.

Soient $E$ et $F$ deux espaces de Hilbert. Un {\em opérateur non borné} est une
application linéaire
$$A : D(A) \to F$$
où $D(A)$ (le domaine de A) est un
sous-espace vectoriel de $E$. En particulier, toute application linéaire
(continue ou non) de $E$ dans $F$ est un opérateur non borné.

Soit $A$ et $B$ deux opérateurs non bornés. On dit que $B$ est un {\em
  prolongement} de $A$, noté $A\subset B$, si $D(A) \subset D(B)$ et
  $B_{|D(A)} = A$.

Un opérateur non borné $A$ est {\em fermé} si son graphe $G(A)= \{(u,A(u)) |
u\in D(A) \}$ est fermé dans $E\times F$, ce qui revient à dire que pour toute
suite $(u_n)$ de $D(A)$ telle que $u_n \to u \in E$ et $A(u_n) \to v \in F$,
on a $u \in D(A)$ et $v=A(u)$.

Si $A$ est à domaine dense dans $E$, on peut définir son {\em adjoint}
$A^* : D(A^*) \subset F \to E$ de la
façon suivante :
$$v \in D(A^*) \Longleftrightarrow \exists w\in
E {\rm \ tel\ que\ } \forall u\in D(A),\ \<u,w\>_E = \<A(u),v\>_F.$$
Comme $D(A)$ est
dense dans $E$, l'élément $w$ (si il existe) est unique ; on pose
$w=A^*(v)$. Remarquons que l'adjoint d'un opérateur est toujours fermé.
On a aussi la propriété évidente (si les opérateurs considérés sont à domaine
dense) $A \subset B \implies B^* \subset A^*$.

Pour définir $A^{**}$, il faut vérifier que $A^*$ est à domaine dense, ce qui
n'est pas toujours le cas. Mais on a la propriété suivante (cf \cite{Riesz} \S
117):

\begin{Prop} Soit $A$ un opérateur non borné de $E$ dans $F$, à domaine
  dense. Alors $A^*$ est à domaine dense si et seulement si $A$ admet un
  prolongement fermé. Dans ce cas, $A^{**}$ est le plus petit prolongement
  fermé de $A$, i.e. si on a $A \subset B$ avec $B$ fermé, alors $A^{**}
  \subset B$.
\end{Prop}

On remarque aussi que le graphe de $A^{**}$ n'est autre que l'adhérence dans
$E\times F$ du graphe de $A$. D'autre part, si $A$ est fermé, on a $A^{**} =
A$. En particulier, dès que cela a un sens, on a toujours $A^{***} = A^*$
(notons au passage que l'on a bien $(A^*)^{**} = (A^{**})^*$).

Si $A$ est injectif, on peut définir son inverse $A^{-1}$ : son domaine n'est
autre que l'image de $A$.

Pour pouvoir définir la composition de deux opérateurs non bornés $A : D(A)
\subset E \to F$ et  $B : D(B) \subset F \to G$, on pose, par définition,
$$D(B\circ A) = \left\{ x \in D(A) | A(x) \in D(B) \right\}.$$
 De même, la somme se définit naturellement sur le domaine
$$D(A+A') = D(A) \cap D(A').$$
Il se peut évidemment que ces domaines soient réduits à l'élément
nul. Cependant, on a le théorème relativement surprenant suivant
(\cite{Riesz}, §118, ou \cite{Rudin}, théorème 13.13) :

\begin{Thm}\label{1+*} Si l'opérateur non borné $A : E \to F$ est fermé et de
  domaine
  dense, alors les opérateurs
$$ B = (A^*\circ A + Id)^{-1},\ C=A\circ (A^*\circ A + Id)^{-1}$$
sont des applications linéaires {\em continues} de $E$ dans $E$ et
de $E$ dans $F$; de plus $||B||\leq 1$, $||C||\leq 1$, et $B$ est
auto-adjointe positive.
\end{Thm}

\bigskip

Maintenant, soit $M$ une variété riemannienne, et soit $E$ et $F$
deux fibrés vectoriels sur $M$, munis de métriques riemanniennes
$(.,.)_E$ et $(.,.)_F.$ On note $C^\infty_0(E)$ (resp.
$C^\infty(E),$ resp. $L^2(E)$) l'espace des sections de E qui sont
$C^\infty$ à support compact (resp. $C^\infty$, resp. $L^2$); de
même pour $F$. La métrique sur $E$ et la forme volume sur $M$ font
de $L^2(E)$ un espace de Hilbert (pour le produit scalaire
${\<f,g\>_E = \int_M (f,g)_E dvol_M}$) dont $C^\infty_0(E)$ est un
sous-espace dense; de même pour $F$.

Soit $A$ un opérateur différentiel agissant sur les sections de
$E$. On le considère comme un opérateur non borné de domaine les
sections $C^\infty$ à support compact, i.e. $$A : C^\infty_0(E)
\to C^\infty_0(F) \subset L^2(F),$$
 et on suppose que $A$ admet un
{\em adjoint formel} $A^t : C^\infty_0(F) \to C^\infty_0(E)$, i.e.
tel que $$\<Au,v\>_F = \<u,A^tv\>_E\ \forall u \in C^\infty_0(E)
{\rm \ et\ } \forall v \in C^\infty_0(F).$$
 On a clairement $A^t \subset A^*$ donc $A^*$ est à domaine dense.

On pose alors $A_{min} = A^{**}$, c'est, on l'a vu,
le plus petit prolongement fermé de $A$. Le graphe de $A^{**}$ est
l'adhérence du graphe de $A$, donc (et on peut prendre ça comme définition)
$$u \in D(A_{min}) \Longleftrightarrow  \exists (u_n) \in C^\infty_0(E) {\rm \
  telle\ que\ } \lim_{n\to \infty} u_n = u {\rm \ et\ que\ la\ suite\ } (Au_n)
  {\rm \ converge\ dans\ } L^2,$$
$A_{min} u$ est alors la valeur de cette limite.

On pose aussi $A_{max} = (A^t)^*$ ; comme $A^t\subset A^*$, on a
$A^{**}\subset A_{max}$ et donc $A \subset A_{max}$. De plus $A^t \subset
(A^t)^{**}=(A_{max})^*$, et, vu la propriété de minimalité de $^{**}$, on en
déduit que $A_{max}$ est le plus grand prolongement
de $A$ dont l'adjoint prolonge aussi $A^t$. Plus précisément,
$$u \in D(A_{max}) \Longleftrightarrow \exists v \in
L^2(F) {\rm \ tel\ que\ } \forall \phi \in C^\infty_0(F),\ \<u,A^t\phi\>_E =
\<v,\phi\>_F,$$
ce qui signifie exactement que $v=Au$ ``au sens des
distributions''. En utilisant des techniques standards d'analyse
(convolution), on montre qu'on peut approcher
$u \in D(A_{max})$ par des sections lisses, i.e. (et on peut prendre ça comme
définition)
$$u \in D(A_{max}) \Longleftrightarrow  \exists (u_n) \in
C^\infty(E) {\rm \ telle\ que\ } \lim_{n \to \infty} u_n = u {\rm \ et\ que\
  la\ suite\ } (Au_n)  {\rm \ converge\ dans\ } L^2$$ ($A_{max} u$ est alors
la valeur de cette limite).

\section{Deux résultats d'intégration par parties sur les
  c\^ones-variétés}\label{secIPP}

Pour faire fonctionner la technique de Bochner nous avons besoin de procéder à
des intégrations par parties. Les deux résultats suivants ainsi que leur
interprétation en termes d'opérateurs non bornés sont à notre disposition.
Le premier théorème d'intégration par parties sur une cône-variété est le
suivant, dû à Cheeger \cite{Cheeger} :

\begin{Thm}\label{Cheeg} Soient $\alpha \in \Omega^p M$ et $\beta \in
  \Omega^{p+1}M$ deux formes $C^\infty$ sur $M$ telles que $\alpha$,
  $d\alpha$, $\beta$, et $\delta\beta$ soient dans $L^2$. Alors
$$\< \alpha, \delta\beta \> = \< d\alpha, \beta\>.$$
\end{Thm}

En fait il faut adapter un tout petit peu la démonstration, ou combiner deux
résultats de l'article cité (cf aussi \cite{HK}, appendice).

Par passage à la limite, il est clair que l'on a encore $$\<
\alpha, \delta_{max}\beta \> = \< d_{max}\alpha, \beta\>$$ quel
que soit $\alpha \in D(d_{max})$ et $\beta \in D(\delta_{max})$.
On en déduit immédiatement (cf aussi \cite{Gaffney}) que :

\begin{Cor}\label{dmax} Les opérateurs $d_{max}$ et $\delta_{max}$ sont
  adjoints l'un de
l'autre; on a $d_{max} = d_{min}$ et $\delta_{max}= \delta_{min}$.
\end{Cor}

\begin{proof} En effet, l'égalité $\< \alpha,\delta_{max}\beta \> = \<
  d_{max}\alpha, \beta\>$ quel que soit $\alpha \in D(d_{max})$ et $\beta \in
  D(\delta_{max})$ implique que $\delta_{max} \subset d_{max}^*$. Or
  $d_{max}^* = \delta_{min} \subset \delta_{max}$. Donc $\delta_{max} =
  \delta_{min} = d_{max}^*$. Le même argument montre que $d_{max} = d_{min} =
  \delta_{max}^*$.
\end{proof}

\bigskip

Le deuxième résultat concerne les tenseurs et non plus les formes
différentielles :

\begin{Thm} \label{ipp} Soient $u \in C^\infty(T^{(r,s)}M)$, $v \in
  C^\infty(T^{(r+1,s)}M)$ tels que $u$, $\n u$, $v$, $\na v$ soient dans
  $L^2$. Alors $$\<u,\na v\> = \<\n u,v\>.$$
\end{Thm}

\begin{proof}
On va démontrer ce résultat en utilisant une méthode similaire à
  celle de Cheeger \cite{Cheeger}. Pour simplifier, nous supposerons que la
  métrique au voisinage de $\Sigma$ est exactement, en coordonnées locales, de
  la forme $g=dr^2+r^2d\theta^2+g_{|\Sigma_i}$, où $\theta$ est définie modulo
  l'angle conique $\alpha$. Le cas général se traite exactement de la même
  façon, les expressions sont juste un peu plus compliquées.

Soit $a$ un réel positif suffisamment petit pour que le
$a$-voisinage fermé de $\Sigma$
  dans $M$ soit tubulaire. Pour $t\leq a$, on pose $U_t$ le $t$-voisinage de
  $\Sigma$ dans $M$, $\Sigma_t = \partial U_t$,
  et $M_t = M\setminus U_t$. Le vecteur $\frac{\partial}{\partial r} = e_r$
  est une normale unitaire en tout point de $\Sigma_t$. Avec ces notations,
  une intégration par parties ( c'est-à-dire la formule de Stokes ) nous donne
  :
$$ \int_{M_t} (g(u,\na v)-g(\n u,v)) = \int_{\Sigma_t} g_{|\Sigma_t}(u,i_{e_r}
v)$$
où $i_{e_r} v = v(e_r,.)$. Le terme de gauche converge vers $\<u,\na v\> - \<\n
u,v\>$ quand $t$ tend vers $0$.
Notons $I_t$ le terme de droite de l'égalité, qui
  correspond au terme de bord. On a alors les inégalités suivantes (la
  notation $|.|$ désigne la valeur absolue aussi bien que la norme ponctuelle
  pour la métrique $g$) :
\begin{eqnarray*}
|I_t| & \leq & \int_{\Sigma_t} |u| |i_{e_r} v| \\
      & \leq & \left(\int_{\Sigma_t} |u|^2\right)^{1/2} \left(\int_{\Sigma_t}
      |i_{e_r} v|^2\right)^{1/2}
\end{eqnarray*}

On va montrer que le fait que $\n u$ soit $L^2$ permet d'avoir une
bonne majoration de $\int_{\Sigma_t} |u|^2$. Et comme $i_{e_r}v$ est $L^2$
(car $v$ l'est aussi), $\int_{\Sigma_t} |i_{e_r}v|^2$ ne
peut pas croître trop vite quand $t$ tend vers $0$. Ce deux résultats nous
permettront d'affirmer que $I_{t_n}$ tend vers $0$ pour une suite $t_n$
tendant vers $0$.

En tout point où $|u| \neq 0$, la fonction $|u|$ est dérivable, et
$d|u|(x) = g(\n _x u, \frac{u}{|u|})$. On pose $$\frac{\partial
|u|}{\partial r} = g(\n_{e_r} u, \frac{u}{|u|})$$ si $|u| \neq 0$,
et $\frac{\partial |u|}{\partial r} =0$ sinon. Il s'agit de la
dérivée partielle distributionnelle de $|u|$, au sens où l'on a,
si les coordonnées autres que $r$ restent fixées, $$|u(t)| -
|u(a)| = \int^t_a \frac{\partial |u|}{\partial r}(r) dr.$$ Or
$|\frac{\partial |u|}{\partial r}| \leq |\n _{e_r} u| \leq |\n
u|$, et donc, si $t\leq a$, $$|u(t)| \leq |u(a)| + \int_t^a |\n u|
dr$$ et $$|u(t)|^2 \leq 2|u(a)|^2 + 2(\int_t^a |\n u| dr)^2.$$

En appliquant l'inégalité de Cauchy-Schwarz il vient
\begin{eqnarray*}
\left(\int_t^a |\n u|dr \right)^2 &\leq &\int_t^a \frac{dr}{r}
\int_t^a r|\n u|^2 dr\\ & \leq & |\ln(\frac{t}{a})| \int_t^a r|\n
u|^2 dr
\end{eqnarray*}

En intégrant sur $\Sigma_t$ on trouve
\begin{eqnarray*}
\int_{\Sigma_t} |u|^2 & \leq & \int_{\Sigma_t} 2|u(a)|^2 +
\int_{\Sigma_t}\left( 2 \ln(\frac{t}{a}) \int_t^a r|\n u|^2
dr\right) \\ & \leq & 2\int_{\Sigma_t} |u(a)|^2 + 2
|\ln(\frac{t}{a})| \int_{\Sigma_t}\int_t^a r|\n u|^2 dr \\ & \leq
& 2\int_{\Sigma}\int_{\theta=0}^\alpha |u(a)|^2td\theta
dvol_\Sigma + 2 |\ln(\frac{t}{a})| \int_{\Sigma}
\int_{\theta=0}^\alpha \left(\int_t^a r|\n u|^2
  dr\right) td\theta dvol_\Sigma\\
& \leq &  2\frac{t}{a}\int_{\Sigma_a} |u|^2 +
2t|\ln(\frac{t}{a})|\int_{\Sigma} \int_{\theta=0}^\alpha \int_t^a
|\n u|^2 rdrd\theta dvol_\Sigma \\ & \leq &
2\frac{t}{a}\int_{\Sigma_a} |u|^2 + 2t|\ln(\frac{t}{a})|
\int_{U_a} |\n u|^2\\ & = & O(t|\ln(t)|)
\end{eqnarray*}

Il reste à contrôler le terme $\int_{\Sigma_t} |i_{e_r}v|^2 \leq
\int_{\Sigma_t} |v|^2$. Comme $v$ est $L^2$, $$\int_0^a
(\int_{\Sigma_t} |v|^2)dt = \int_{U_a} |v|^2 < +\infty,$$ et donc
la fonction $t \mapsto \int_{\Sigma_t} |v|^2$ est intégrable sur
$]0,a]$. Or la fonction $(t\ln(t))^{-1}$ n'est pas intégrable en
$0$. On en déduit donc qu'il existe une suite $t_n$ tendant vers
$0$ pour laquelle $$\int_{\Sigma_{t_n}} |v|^2 =
o((t_n\ln(t_n))^{-1}).$$ En combinant avec la majoration obtenue
pour $\int_{\Sigma_t} |u|^2$, on en déduit immédiatement que
$${\lim_{n\to +\infty} I_{t_n} = 0}.$$ Or $I_t \to \<u,\na v\> -
\<\n u,v\>$ quand $t \to 0$; on a donc bien  $\<u,\na v\> = \<\n
u,v\>.$
\end{proof}

\bigskip

\begin{Cor} \label{ipp1} On considère $\n$ comme un opérateur non borné $\n :
  C^\infty_0(T^{(r,s)}M) \to C^\infty_0(T^{(r+1,s)}M)$. Alors $\forall u \in
  D(\n_{max}),$ $\forall v \in D(\na),$ on a l'égalité
$$\<\n_{max} u,v\> = \<u, \na v\>.$$ Ceci implique en particulier que
 $\n_{min} = \n_{max}$ et que $\n^t_{min} = \n^t_{max} = \na$.
\end{Cor}

\begin{proof}
La première égalité se démontre directement en prenant des suites
régularisantes pour $u$ et pour $v$ : en effet on a vu dans la section
précédente que si $u \in D(\n_{max})$ (resp. $v \in D(\na)$), il existe une
suite $u_n \in C^\infty(T^{(r,s)}M)$ (resp. $v_n \in C^\infty(T^{(r+1,s)}M)$)
telle que $\lim_{n\to\infty} u_n = u$ et $\lim_{n\to\infty} \n u_n = \n_{max}
u$ dans $L^2$ (resp. $\lim_{n\to\infty} v_n = v$ et $\lim_{n\to\infty} \na v_n
= \na v$). Alors
$$\<u,\na v\> = \lim \<u_n,\na v_n\> = \lim \<\n u_n, v_n\> = \< \n_{max} u,
v\>.$$
La suite se démontre comme le corollaire \ref{dmax}.
\end{proof}

Dans la suite les notations $\n$ et $\na$ désigneront donc (sauf exception)
les opérateurs $\n_{max} ( = \n_{min} )$ et $\n^t_{max} ( = \n^t_{min} = \na).$

On employera fréquemment le corollaire suivant, simple reformulation du
précédent :

\begin{Cor}\label{approx} Soit $u$ appartenant à $D(\n)$, c'est-à-dire tel que
  $u$ et $\n u$
  sont $L^2$. Alors il existe une suite $(u_n)$, $C^\infty$ à support
  compact, telle que $u_n \to u$ et $\n u_n \to \n u$ dans $L^2$ quand $n \to
  \infty$.
\end{Cor}

\begin{proof} C'est juste la définition de $u \in D(\n_{min})$.
\end{proof}

\bigskip

{\em Remarque :} dans les théorèmes ci-dessus, la condition $L^2$ paraît
naturelle, ne serait-ce que pour s'assurer de l'existence des termes du type
$\<\n u,v\>$. Cependant il est intéressant de noter que la démonstration
donnée du théorème \ref{ipp} ne fonctionne pas avec des hypothèses du type
$u$, $\n u \in L^p$ et $v$, $\na v \in L^q$, avec $p$ et $q$ des exposants
conjugués. La condition $L^2$ est donc plus importante qu'il n'y paraît.

\section{Etude de l'équation de normalisation}\label{secLap}

Dans toute cette section, nous supposerons que la métrique conique $g$ est {\em
  hyperbolique}.

Comme on l'a vu dans la section \ref{secEin}, pour se débarrasser des
déformations triviales on cherche à imposer la condition de jauge de
Bianchi. Montrer qu'une déformation infinitésimale peut se mettre sous une
forme normalisée vérifiant la condition de jauge revient à résoudre l'équation
de normalisation :
$$\beta \circ \delta^* \alpha = \beta h_0.$$

Cette équation peut se mettre sous une forme plus lisible. Pour cela, on
utilise le fait que $$\n \alpha = \delta^* \alpha + \frac{1}{2}d \alpha$$ (il
s'agit juste de la décomposition du $2$-tenseur $\n \alpha$ en partie
symétrique et anti-symétrique), que $\delta$ est toujours la restriction de
$\na$ au sous-fibré correspondant, et donc que $$\delta \alpha = \na \alpha =
- \tr \n \alpha = - \tr \delta^* \alpha,$$
la trace de $d \alpha$ étant nulle puisque $d\alpha$ est anti-symétrique. On
 obtient alors
\begin{eqnarray*}
2 \beta (\delta^* \alpha) & = & 2 \delta \delta^* \alpha + d \tr \delta^*
\alpha \\
& = & 2 \na (\n \alpha - \frac{1}{2}d \alpha) - d \delta \alpha \\
& = & 2 \na \n \alpha - \delta d \alpha - d \delta \alpha \\
& = & 2 \na \n \alpha - \Delta \alpha.
\end{eqnarray*}
 Ici $\Delta = d\delta + \delta d$ est l'opérateur de Laplace-Beltrami sur les
 1-formes. Or $\Delta$ et $\na \n$ (parfois nommé {\em laplacien de connexion}) sont
 reliés par la classique formule de Weitzenböck $$\Delta\alpha = \na \n \alpha
 + ric(\alpha),$$
cf \cite{Besse} \S 1.155. En utilisant cette formule et le fait que la
 métrique $g$ est hyperbolique, on trouve
\begin{eqnarray*}
2 \beta (\delta^* \alpha) & = & \Delta \alpha + 2(n-1) \alpha \\
& = & \na \n \alpha + (n-1) \alpha.
\end{eqnarray*}

On est donc amener à étudier l'opérateur $L : \alpha  \mapsto \na \n \alpha +
(n-1) \alpha$.

\subsection{Premières propriétés}\label{LapDeb}

La première chose à remarquer sur $L$ est qu'il est {\em elliptique}. En
particulier, si $\phi$ est $C^\infty$ et que $L\alpha = \phi$ au sens des
distributions, alors $\alpha$ est $C^\infty$. Malheureusement, le caractère
singulier d'une cône-variété nous empêche d'utiliser directement les
inégalités de type Schauder ou G\r{a}rding. Par exemple, on peut montrer qu'il
existe des 1-formes $\alpha$ appartenant à $L^2$
telles que $L\alpha = 0$ au sens des distributions avec $\n \alpha$ qui n'est
pas dans $L^2$.

\bigskip

Il est clair que $L$, vu comme un opérateur non borné $C^\infty_0(T^*M) \to
C^\infty_0(T^*M)$, est formellement symétrique : avec les notations de la
section \ref{rapop}, $L^t = L$. Malheureusement, il est possible de montrer
que dès que la dimension de
notre cône-variété est supérieure à 2, l'opérateur $L$ n'est pas
essentiellement auto-adjoint, i.e. $L_{min} \neq L_{max}$ (ou si l'on préfère,
$L^{**} \neq L^*$). On va donc étudier
des extensions auto-adjointes $\bar{L}$ de $L$, avec $L_{min} \subset \bar{L}
\subset L_{max}$.

Le théorème \ref{1+*} nous donne une première telle extension : toujours avec
les
conventions de la section \ref{rapop}, l'opérateur $\na \circ \n + (n-1) Id$
est auto-adjoint et inversible. Son domaine est par définition
$$D = \left\{ \alpha \in D(\n) |\ \n \alpha \in D(\na)\right\} = \left\{ \alpha
  \in  L^2 |\ \n \alpha, \ \na \n \alpha \in L^2 \right\}$$
(dans la deuxième définition, il faut considérer $\n$ et $\na \n$ au sens des
  distributions).

On va maintenant introduire une deuxième extension auto-adjointe. On doit à
Gaffney \cite{Gaffney} le résultat général suivant : les opérateurs $d_{min}
\delta_{max} +
\delta_{max} d_{min}$ et $d_{max} \delta_{min} + \delta_{min} d_{max}$
(encore avec les conventions de la section \ref{rapop}) sont toujours
auto-adjoints. Or d'après le corollaire \ref{dmax}, on a $d_{min} = d_{max}$
et $\delta_{min} = \delta_{max}$; les deux opérateurs ci-dessus sont donc les
mêmes sur une cône-variété. On en déduit que l'opérateur $d_{max} \delta_{max}
+ \delta_{max} d_{max} + 2(n-1) Id$, défini sur le domaine
\begin{eqnarray*}
D'& = & \left\{\alpha \in D(d_{max}) \cap D(\delta_{max}) |\ d_{max} \alpha
  \in D(\delta_{max}) \mathrm{\ et\ } \delta_{max} \alpha \in D(d_{max})
  \right\} \\
 & = & \left\{ \alpha \in  L^2 |\ d \alpha,\ \delta \alpha,\ d\delta\alpha,\
  \delta d\alpha \in L^2 \right\},
\end{eqnarray*}
est positif auto-adjoint et donc inversible (encore une fois dans la deuxième
définition, il faut considérer les opérateurs au sens des distributions). Nous
montrerons plus
loin que ces deux domaines $D$ et $D'$ sont en fait
confondus quand tous les angles coniques sont inférieurs à $2\pi$.

\subsection{Expression du laplacien de connexion en coordonnées
  cylindriques}\label{vois1}

On va maintenant sauter à pieds joints dans les calculs.
Soit $a$ un réel positif suffisament petit pour que le $a$-voisinage fermé de
$\Sigma$ dans $M$ soit un voisinage tubulaire. Si $r$ est plus petit que $a$,
on note $U_r$ le $r$-voisinage de $\Sigma$ dans $M$ et $\Sigma_r$ le bord de
$U_r$.

Par définition, si $x$ est un point de $\Sigma$, il existe un voisinage $V$ de
$x$ dans $U_a$ et un voisinage $U$ de $x$ dans $\Sigma$, tels que $U=V\cap
\Sigma$ et $V \simeq U\times D^2$, et dans les coordonnées cylindriques locales
adaptées à la décomposition $V \simeq U\times D^2$, la métrique est de la
forme
$$g=dr^2 + \sinh(r)^2 d\theta^2 + \cosh(r)^2 g_\Sigma,$$
où $\theta$ est défini non pas modulo $2\pi$ mais modulo l'angle conique
$\alpha$. On utilisera les notations
suivantes : $e^r= dr$, $e^\theta = \sinh(r) d\theta$, $e_r = (e^r)^\sharp =
\frac{\d}{\d r}$, et $e_\theta = (e^\theta)^\sharp =
\frac{1}{\sinh(r)}\frac{\d}{\d \theta}$.

Soit $u$ une section de $T^*M$. Au voisinage de $\Sigma$ on peut faire une
décomposition orthogonale et on écrit
$$u = f e^r + g e^\theta + \omega,$$
 avec $f$, $g$, deux fonctions de $M$ dans $\mathbb{R}$ (ou $\mathbb{C}$, on
 sera souvent amené dans la suite à complexifier les fibrés sur lesquels on
travaille), et $\omega$ une 1-forme. On
remarque que bien que les coordonnées ne soient que
locales, les formes $e^r$ et $e^\theta$ sont bien définies sur tout $U_a$,
ainsi que la décomposition orthogonale précédente.

Au vu de la forme de notre voisinage tubulaire, sur tout ouvert $V$ de $U_a$ du
type ci-dessus et suffisament petit, on peut définir localement des
champs de vecteur $e_1,\ldots e_{n-2}$ de telle sorte que
$(e_r,e_\theta,e_1,\ldots e_{n-2})$ forme un repère mobile orthonormé (local),
vérifiant
$$\n_{e_r} e_k = \n_{e_\theta} e_k = 0$$
pour tout $k$ dans $1\ldots n-2$. On définit de même des 1-formes locales
$e^1,\ldots e^{n-2}$ telles
que \\ $(e^r,e^\theta,e^1,\ldots e^{n-2})$ soit le repère mobile dual du
précédent.

Avant de commencer les calculs, introduisons encore quelques notations. On
note $N$ le (sous-)fibré vectoriel au-dessus de $U_a$, dont la fibre
au-dessus de $x \in U_a$ est le sous-espace vectoriel de $T^*_x M$
orthogonal à $e^\theta$ et $e^r$, et $N^*$ le (sous-)fibré vectoriel au-dessus
de $U_a$, dont la fibre au-dessus de $x \in U_a$ est le sous-espace vectoriel
de $T_x M$ orthogonal à $e_\theta$ et $e_r$. La $1$-forme $\omega$ introduite
plus haut est naturellement une section de $N$. Les sections $(e_1,\ldots,
e_{n-2})$ forment localement une base de $N^*$, de même pour $(e^1,\ldots,
e^{n-2})$ et $N$. Si $s$ est section de $N^*$,
et $t$ une section de $N$ ou de $N^*$, on note $(\n_\Sigma)(s,t)$ la
projection orthogonale sur $N$ ou sur $N^*$ de $\n_s t$.

Si $f$ est une fonction de $U_a$, on pose $$d_\Sigma f =
\sum_{k=1}^{n-2} (e_k.f)e^k,$$ et $$\Delta_\Sigma f =  \cosh^2 r
\sum_{k=1}^{n-2} (\n_\Sigma)(e_k,e_k).f - e_k.e_k.f$$ (c'est à un
facteur près l'opposé de la trace de la hessienne de $f$
restreinte à $N^*$). Ces deux opérateurs sont indépendants du
choix des $e_k$. En fait, avec les notations ci-dessus, dans $V$
on a une identification, à $r$ et $\theta$ fixé, de $U\times
\{r,\theta\}$ à $U\subset \Sigma$, et $N^*$ et $N$ restreints à
$U\times \{r,\theta\}$ s'identifient de la même façon à $TU$ et
$T^*U$. Les opérateurs ci-dessus correspondent via ces
identifications à la différentielle et au laplacien de $U\subset
\Sigma$.

Il en est de même pour $\n_\Sigma$, et pour les deux opérateurs suivants. Si
$\omega$ est une section de $N$, on pose
\begin{eqnarray*}\delta_\Sigma (\omega) & = & - \cosh^2 r
\sum_{k=1}^{n-2} \n_\Sigma(e_k,\omega)(e_k) \\
& = & \cosh^2 r \sum_{k=1}^{n-2}
  \omega(\n_\Sigma(e_k,e_k)) - e_k.\omega(e_k),
\end{eqnarray*}
et
$$(\na \n)_\Sigma \omega = \cosh^2 r \sum_{k=1}^{n-2}
\n_\Sigma(\n_\Sigma(e_k,e_k), \omega) - \n_\Sigma (e_k, \n_\Sigma( e_k,
\omega)),$$
qui correspondent à la codifférentielle et au laplacien de connexion pour les
1-formes de $\Sigma$.

\bigskip

On est maintenant armé pour le calcul explicite de $\na \n u$. En utilisant
notre repère mobile, on a
\begin{eqnarray*}
\na \n u & = & -\n_{e_r}\n_{e_r}u - \n_{e_\theta}\n_{e_\theta}u +
\n_{\n_{e_r}e_r} u + \n_{\n_{e_\theta}e_\theta}u + \sum_{k=1}^{n-2}
\n_{\n_{e_k}e_k} u -\n_{e_k}\n_{e_k}u.
\end{eqnarray*}
Comme la métrique conique est hyperbolique, on a les expressions suivantes :
$$\n_{e_r}e_r = 0,\ \n_{e_\theta}e_\theta = -\frac{1}{\tanh(r)}e_r, {\rm \ et\
  } \n_{e_k}e_k = \n_\Sigma(e_k,e_k) - \tanh(r)\,e_r.$$
On trouve alors que
$$\na \n u  =  -\n_{e_r}\n_{e_r}u - \n_{e_\theta}\n_{e_\theta}u -
 (\frac{1}{\tanh(r)} + (n-2)\tanh(r))\n_{e_r} u + \sum_{k=1}^{n-2}
\n_{\n_\Sigma(e_k,e_k)} u -\n_{e_k}\n_{e_k}u.$$

En remplaçant $u$ par $f e^r + g e^\theta + \omega$, un calcul explicite nous
donne
l'expression suivante pour les composantes de
$\na \n u$; selon $e^r$ :

$$ -\frac{\partial^2 f}{\partial r^2} -
\frac{1}{\sinh(r)^2}\frac{\partial^2 f}{\partial \theta ^2}
- \left(\frac{1}{\tanh(r)} + (n-2) \tanh(r)\right) \frac{\partial f}{\partial
  r} + \left(\frac{1}{\tanh(r)^2} + (n-2) \tanh(r)^2\right)f \coupeq
 + \frac{2}{\sinh(r) \tanh (r)} \frac{\partial g}{\partial \theta} +
 \frac{1}{\cosh(r)^2} \Delta_\Sigma f
+ \frac{2 \tanh(r)}{\cosh(r)^2} \delta_\Sigma \omega, $$

selon $e^\theta$ :

$$ -\frac{\partial^2 g}{\partial r ^2} -
\frac{1}{\sinh(r)^2}\frac{\partial^2 g}{\partial \theta ^2}
-\left(\frac{1}{\tanh(r)} + (n-2) \tanh(r)\right)\frac{\partial g}{\partial r}
+ \frac{g}{\tanh(r)^2} - \frac{2}{\sinh(r) \tanh
(r)}\frac{\partial f}{\partial \theta} + \frac{1}{\cosh(r)^2} \Delta_\Sigma
g,$$

et selon la composante incluse dans $N$ :

$$ -\n_{e_r} \n_{e_r} \omega - \n_{e_\theta}\n_{e_\theta} \omega
 - \left(\frac{1}{\tanh(r)} + (n-2) \tanh(r)\right) \n_{e_r} \omega +
 \tanh(r)^2 \omega \coupeq
 + \frac{1}{\cosh(r)^2} (\na \n)_\Sigma \omega - 2\tanh(r) \,d_\Sigma f$$

\bigskip

Pour pouvoir manipuler cette expression, nous allons effectuer dans la section
suivante une sorte de décomposition en séries de Fourier généralisées.

\subsection{Décomposition en série de Fourier généralisée}\label{vois2}

On sait qu'au voisinage du lieu singulier, la métrique $g$ se met
localement sous la forme $$g = dr^2 + \sinh(r)^2 d\theta^2 +
\cosh(r)^2 g_\Sigma.$$ Si la coordonnée $\theta$ était définie
(toujours modulo l'angle conique $\alpha$) sur tout un voisinage
du lieu singulier, on pourrait faire des décompositions en séries
de Fourier, du type $$f(r,\theta,z) = \sum f_n(r,z) \exp(2i\pi
n\theta/\alpha).$$ Mais en général la coordonnée d'angle $\theta$
n'est définie que localement, ce qui empêche d'écrire de telles
décompositions. On va donc procéder à une autre sorte de
décomposition; on obtiendra finalement des écritures du type
$$f(r,\theta,z) = \sum f_n(r) \psi_n(\theta,z)$$
 où les $(\psi_n)$
forment une base hilbertienne bien choisie du bord d'un voisinage
tubulaire du lieu singulier.

\subsubsection{Une base hilbertienne adaptée}

Pour faciliter les calculs, nous serons amener à complexifier les fibrés
usuels. Comme précédemment, on choisit un réel positif $a$ suffisament petit
pour que le $a$-voisinage fermé de $\Sigma$ dans $M$ soit un voisinage
tubulaire. Si $r$ est inférieur ou égal à $a$,
on note $U_r$ le $r$-voisinage de $\Sigma$ dans $M$ et $\Sigma_r$ le bord de
$U_r$.

On va particulièrement s'intéresser à la sous-variété $\Sigma_a$. Tout point
$x$ de $\Sigma_a$ admet un voisinage $\mathcal{V}$ de la forme $U \times S^1$,
où $U$ est un ouvert de
$\Sigma$. Dans ce voisinage, la métrique de $\Sigma_a$, induite par celle
de $M$, s'exprime comme une métrique produit; plus précisement on a, dans
les coordonnées adaptées,
$$g_a = \sinh(a)^2 d\theta^2 + \cosh(a)^2 g_\Sigma.$$
Dans cette sous-section, et seulement dans celle-ci, les notations $\n$,
$\na$, $\Delta$, etc. désigneront les opérateurs correspondants pour la
métrique $g_a$.

Pour pouvoir faire les décompositions voulues, on veut trouver une
``bonne'' base hilbertienne sur $\Sigma_a$, pour les fonctions
comme pour les $1$-formes, ou plus précisément pour les sections
du sous-fibré $N$ défini précédemment. On se propose donc de
démontrer le résultat suivant :

\bigskip

\begin{Prop}\label{bashib} Il existe une base hilbertienne $(\psi_j)_{j\in
    \mathbb{N}}$ du complexifié de $L^2(\Sigma_a)$, telle que pour tout indice $j$, il
  existe un réel $\lambda_j$ et un entier relatif $p_j$, vérifiant
$$\lambda_j \geq \frac{p_j^2 \beta^2}{\sinh(a)^2},$$
 pour lesquels
$$\begin{cases} \Delta \psi_j = \lambda_j \psi_j \\
e_\theta.\psi_j = \frac{ip_j\beta}{\sinh(a)}\psi_j.
\end{cases}$$

Soit $J$ l'ensemble des $j$ pour lesquels $\lambda_j > \frac{p_j^2
  \beta^2}{\sinh(a)^2}$. Il existe une base hilbertienne $(\phi_j)_{j\in J}
  \cup (\varphi_j)_{j\in \mathbb{N}}$ du complexifié de $L^2(N)$, telle que :
\begin{itemize}
\item pour tout indice $j$ appartenant à $J$, $\phi_j = \left(\lambda_j - \frac{p_j^2
      \beta^2}{\sinh(a)^2}\right)^{-1/2} d_\Sigma \psi_j$, et donc
$$\begin{cases} \na \n \phi_j = \left(\lambda_j +
    \frac{n-3}{\cosh(a)^2}\right)\phi_j \\
\n_{e_\theta}\phi_j = \frac{ip_j\beta}{\sinh(a)}\phi_j\\
\delta_\Sigma \phi_j = \cosh(a)^2 \left(\lambda_j - \frac{p_j^2
      \beta^2}{\sinh(a)^2}\right)^{1/2} \psi_j;
\end{cases}$$
\item pour tout indice $j \in \mathbb{N}$, il existe un réel $\mu_j$ et un
  entier relatif $p'_j$, pour lesquels
$$\begin{cases} \na \n \varphi_j = \mu_j \varphi_j \\
\n_{e_\theta}\varphi_j = \frac{ip'_j\beta}{\sinh(a)}\varphi_j,
\end{cases}$$
et on a de plus $\delta_\Sigma \varphi_j = 0$.
\end{itemize}
\end{Prop}

\medskip

\begin{proof}
Soit $(e_1,\ldots,e_{n-2})$ un repère
mobile orthonormé de $U \simeq U \times \{\theta\}$, pour la métrique
$\cosh(a)^2 g_\Sigma$. C'est la restriction à $\Sigma_a$ du repère mobile
local défini dans la section précédente. Alors $(e_\theta,e_1,\ldots,e_{n-2})$
est un repère
mobile orthonormé de $\mathcal{V}$, vérifiant
$$\n_{e_\theta} e_k = \n_{e_k} e_\theta = 0$$
pour $k=1\ldots n-2$ (rappelons qu'ici et dans la suite de la preuve, on
réemploie pour simplifier la
notation $\n$ pour la connexion de
Levi-Cività pour la métrique $g_a$ de $\Sigma_a$, à ne pas confondre avec la
connexion de Levi-Cività de la métrique $g$).

Intéressons-nous maintenant au laplacien $\Delta$ sur $(\Sigma_a,g_a)$. La
sous-variété $\Sigma_a$ étant compacte, on peut utiliser le théorème de
décomposition spectrale des opérateurs elliptiques auto-adjoints pour montrer
qu'il
existe une base hilbertienne de $L^2(\Sigma_a)$, formée de fonctions propres
du laplacien. De plus chaque sous-espace propre est de dimension finie,
et chaque fonction propre est $C^\infty$.

Si $f$ est une fonction sur $\Sigma_a$, en utilisant le repère
mobile ci-dessus, on obtient $$\Delta f = - e_\theta.e_\theta.f +
\sum_{k=1}^{n-2} \n_{e_k}e_k.f - e_k.e_k.f\,.$$

On peut vérifier sur cette expression que $e_\theta.\Delta f = \Delta
(e_\theta.f)$, donc que $\frac{\d}{\d \theta}(\Delta f) = \Delta (\frac{\d}{\d
  \theta} f)$, c'est-à-dire que $\Delta$ et $\frac{\d}{\d \theta}$ commutent
(rappelons que le champ de vecteur $e_\theta$ est bien défini dans
tou $U_a$, ainsi donc que $\frac{\d}{\d \theta}$). De plus une
intégration par parties évidente donne, si $f$ et $g$ sont deux
fonctions $C^1$ sur $\Sigma_a$, $$\int_{\Sigma_a} \frac{\d f}{\d
\theta}g = - \int_{\Sigma_a} f\frac{\d g}{\d
  \theta}.$$

Soit donc $\lambda$ une valeur propre du laplacien et $E_\lambda$ le
sous-espace propre associé. On a vu que $E_\lambda$ était de dimension fini et
composé de fonctions $C^\infty$. Comme $\Delta$ et $\frac{\partial}{\partial
  \theta}$ commutent, la restriction de ce dernier à $E_\lambda$ est un
endomorphisme de $E_\lambda$; et d'après ce qui précède cet endomorphisme est
anti-symétrique pour le produit scalaire $L^2$. Comme $E_\lambda$ est
de dimension finie, on peut donc (en passant dans les complexes) trouver une
base orthonormée $(\psi_j)$ de $E_\lambda$, formée de fonctions propres de
$\frac{\partial}{\partial \theta}$, c'est-à-dire que l'on a :
$$\begin{cases}
  \Delta \psi_j = \lambda \psi_j \\ \frac{\partial}{\partial \theta}\psi_j =
  i \mu_j \psi_j.
\end{cases}$$
Cette dernière équation implique que, dans les coordonnées
$(z,\theta)$ adaptées à $\mathcal{V} \simeq U \times S^1$,
$$\psi_j(z,\theta) = \exp(i\mu_j\theta)\phi_j(z),$$ où $\phi_j$
est une fonction ne dépendant pas de la variable $\theta$. Comme
$\theta$ est définie modulo l'angle conique $\alpha$, on en déduit
que $\mu_j \in \frac{2\pi}{\alpha}\mathbb{Z}$, c'est-à-dire qu'il
existe un entier $p_j$ tel que $$\psi_j(z,\theta) =
\exp(\frac{2i\pi p_j}{\alpha}\theta)\phi_j(z),$$
 et donc en particulier
$$\frac{\partial}{\partial \theta}\psi_j = \frac{2i\pi
p_j}{\alpha} \psi_j.$$

Si on exprime à nouveau le laplacien à l'aide du repère mobile, on obtient
\begin{eqnarray*} \Delta \psi_j & = & \lambda \psi_j \\
  & = &  - e_\theta.e_\theta.\psi_j +
\sum_{k=1}^{n-2} \n_{e_k}e_k.\psi_j - e_k.e_k.\psi_j \\
& = & -\frac{1}{\sinh(a)^2}\frac{\d^2}{\d \theta ^2}\psi_j + \sum_{k=1}^{n-2}
\n_{e_k}e_k.\psi_j - e_k.e_k.\psi_j \\
& = & \frac{p_j^2 \beta^2}{\sinh(a)^2} \psi_j +
\sum_{k=1}^{n-2} \n_{e_k}e_k.\psi_j - e_k.e_k.\psi_j
\end{eqnarray*}
avec la notation $\beta = \frac{2\pi}{\alpha}$. Désignons, comme
dans la section précédente, par $\Delta_\Sigma$ l'opérateur $$f
\mapsto \cosh(a)^2 \sum_{k=1}^{n-2} (\n_{e_k}e_k.f - e_k.e_k.f);$$
si $f$ est une fonction définie sur $U \times \{\theta\}$, alors
$\Delta_\Sigma f$ est le laplacien de $f$ pour la métrique obtenue
en identifiant $U \times \{\theta\}$ à $U \subset \Sigma$. On a
alors $$\Delta_\Sigma (\psi_j) = \cosh(a)^2 \left(\lambda - \frac{p_j^2
  \beta^2}{\sinh(a)^2}\right) \psi_j.$$

\bigskip

Passons maintenant aux 1-formes sur $\Sigma_a$. Rappelons qu'ici
les notations $\n$ et $\na$ désignent la connexion de Levi-Cività
et son adjoint pour la métrique $g_a$ de $\Sigma_a$. Le laplacien
de connexion $\na \n$ s'exprime alors à l'aide du repère mobile de
la façon suivante : $$\na \n \eta = - \n_{e_\theta} \n_{e_\theta}
\eta + \sum_{k=1}^{n-2} \n_{\n_{e_k}e_k} \eta - \n e_k \n e_k
\eta.$$ Si on décompose $\eta$ orthogonalement, $\eta = fe^\theta
+ \omega$, alors $$\na \n \eta = (\Delta f) e^\theta + \na \n
\omega,$$ et cette décomposition est à nouveau orthogonale. Dit
autrement, si $\omega$ est une 1-forme sur $\Sigma_a$, en tout
point perpendiculaire à $e^\theta$, alors $\na \n \omega$ est
aussi en tout point perpendiculaire à $e^\theta$. On va donc
considérer le sous-fibré vectoriel $N \subset T^*\Sigma_a$, dont
la fibre au-dessus de $x \in \Sigma_a$ est le sous-espace
vectoriel de $T^*_x \Sigma_a$ orthogonal à $e^\theta$; c'est la
restriction à $\Sigma_a$ du fibré $N$ défini à la section
précédente. L'opérateur $\na \n$ se restreint ainsi à un opérateur
non borné de $L^2(N)$ dans lui-même.

Si $\omega$ est une section de $N$, alors $\n_{e_\theta} \omega$ est encore
une section de $N$. Comme pour les fonctions, on vérifie sur
l'expression ci-dessus que $\n_{e_\theta} (\na \n \omega) = \na \n
(\n_{e_\theta} \omega)$, c'est-à-dire que $\na \n$ et $\n_{e_\theta}$
commutent. L'opérateur $\n_{e_\theta}$ est à nouveau anti-symétrique pour le
produit scalaire $L^2$ : en effet, pour $\phi$, $\phi'$ deux sections de $N$,
on a
\begin{eqnarray*} \int_{\Sigma_a} e_\theta.g_a(\phi,\phi') & = & 0 \\
& = & \int_{\Sigma_a} g_a(\n_{e_\theta}\phi,\phi ') +\int_{\Sigma_a}
g_a(\phi,\n_{e_\theta} \phi')
\end{eqnarray*}
En utilisant ces deux propriétes et le théorème de décomposition
spectrale des opérateurs elliptiques auto-adjoints, on montre,
comme dans le cas des fonctions, qu'il existe un base hilbertienne
$(\phi_j)_{j\in\mathbb{N}}$ de $L^2(N)$, telle que pour tout $j$,
on ait : $$\begin{cases} \na \n \phi_j = \lambda_j \phi_j \\
\n_{\frac{\d}{\d \theta}}
  \phi_j = i \mu_j \phi_j. \end{cases}$$

Pour $k=1\ldots n-2$, on note comme précédemment $e^k$ la forme
duale de $e_k$. Les sections locales $e^1,\ldots,e^{n-2}$ forment
alors une base de $N$ sur $\mathcal{V}$. En décomposant dans cette
base $$\phi_j = \sum_{k=1}^{n-2} a_k(z,\theta)e^k,$$ la dernière
équation ci-dessus donne $$\sum_{k=1}^{n-2} \frac{\d}{\d
\theta}(a_k(z,\theta))e^k = i \mu_j \sum_{k=1}^{n-2}
a_k(z,\theta)e^k,$$ qui s'intègre en $$a_k(z,\theta) =
b_k(z)\exp(i\mu_j\theta)$$ pour $k=1\ldots n-2$. La coordonnée
$\theta$ étant définie modulo l'angle conique $\alpha$, on en
déduit encore une fois que $\mu_j \in
\frac{2\pi}{\alpha}\mathbb{Z}$, c'est-à-dire qu'il existe un
entier $p_j$ tel que $\n_{\frac{\d}{\d \theta}} \phi_j = i p_j
\beta \phi_j$ (avec toujours $\beta = \frac{2\pi}{\alpha}$).

\bigskip

Rappelons maintenant quelques notations. Si $f$ est une fonction sur
$\Sigma_a$, on a
$$d_\Sigma f = \sum_{k=1}^{n-2} (e_k.f) e^k,$$
 c'est une section de $N$. Il s'agit en chaque point $(x,\theta) \in U\times
 S^1 \simeq \mathcal{V} \subset \Sigma_a$ de la différentielle de la
 restriction de
 $f$ à $U\times\{\theta\}$. Ensuite, si $\omega$ est une section de $N$, on
 a
\begin{eqnarray*}\delta_\Sigma (\omega) & = & - \cosh(a)^2
\sum_{k=1}^{n-2} (\n_{e_k}u)(e_k) \\ & = & \cosh(a)^2 \sum_{k=1}^{n-2}
  u(\n_{e_k}e_k) - e_k.u(e_k),\end{eqnarray*} et
$$(\na \n)_\Sigma \omega = \cosh(a)^2 \sum_{k=1}^{n-2} \n_{\n_{e_k}e_k} \omega
- \n_{e_k} \n_{e_k} \omega.$$
Si l'on se restreint à $U\times \{\theta\}$, il s'agit de la codifférentielle
et du laplacien de connexion sur les 1-formes, pour la métrique obtenue en
identifiant $U\times \{\theta\}$ à $U \subset \Sigma$.

Soit $f$ une fonction sur $\Sigma_a$. Avec ces notations, on a :
$$ \na \n (d_\Sigma f) = - \n_{e_\theta} \n_{e_\theta}
  (d_\Sigma f) + \frac{1}{\cosh(a)^2} (\na \n)_\Sigma (d_\Sigma
  f).$$
Or \begin{eqnarray*} \n_{e_\theta} \n_{e_\theta}
  (d_\Sigma f) & = & \n_{e_\theta} \n_{e_\theta} \sum_{k=1}^{n-2}
  (e_k.f) e^k \\
& = & \sum_{k=1}^{n-2} (e_\theta.e_\theta.(e_k.f)) e^k \\
& = & \sum_{k=1}^{n-2} e_k.(e_\theta.e_\theta.f) e^k \\
& = & d_\Sigma (e_\theta.e_\theta.f).
\end{eqnarray*}
De plus, la métrique de $U\subset \Sigma$ est hyperbolique, et on peut donc
utiliser la formule de Weitzenböck suivante, valable pour les 1-formes :
$$(\na \n)_\Sigma = \Delta_\Sigma + (n-3)Id = d_\Sigma \delta_\Sigma +
\delta_\Sigma d_\Sigma + (n-3)Id.$$
C'est la même formule qu'au début de cette section, cf \cite{Besse} \S 1.155.
On en déduit que $$(\na \n)_\Sigma \circ d_\Sigma = (\Delta_\Sigma +
(n-3)Id)\circ d_\Sigma = d_\Sigma \circ (\Delta_\Sigma + (n-3)Id).$$
Par suite,
$$(\na \n)_\Sigma (d_\Sigma f) = d_\Sigma(\Delta_\Sigma(f)) +
(n-3)d_\Sigma(f),$$
et donc
\begin{eqnarray*}\na \n (d_\Sigma f) & = & - d_\Sigma (e_\theta.e_\theta.f) +
\frac{1}{\cosh(a)^2}\left(d_\Sigma(\Delta_\Sigma(f)) +
  (n-3)d_\Sigma(f)\right)\\
 & = & d_\Sigma(\Delta f) + \frac{n-3}{\cosh(a)^2} d_\Sigma f.
\end{eqnarray*}
On obtient finalement, si $\psi_j$ est une des fonctions définies
plus haut (c'est-à-dire telle que  $\Delta \psi_j = \lambda
\psi_j$ et $\frac{\d}{\d \theta} \psi_j = i p_j \beta \psi_j$) :
\begin{eqnarray*} \na \n (d_\Sigma \psi_j) &=& d_\Sigma (\lambda \psi_j) +
  \frac{n-3}{\cosh(a)^2} d_\Sigma \psi_j\\
&=& \left(\lambda+\frac{n-3}{\cosh(a)^2}\right)d_\Sigma \psi_j.
\end{eqnarray*}
Par conséquent $d_\Sigma \psi_j$ est un vecteur propre de $\na \n$. De plus
$\n_{\frac{\d}{\d \theta}} (d_\Sigma \psi_j) = d_\Sigma (\frac{\d}{\d
  \theta}\psi_j) = i p_j \beta \, d_\Sigma \psi_j$.

\bigskip

De la même manière, si $\omega$ est une section de $N$, alors
$$\Delta (\delta_\Sigma \omega) = - e_\theta.e_\theta.(\delta_\Sigma \omega) +
\frac{1}{\cosh(a)^2} \Delta_\Sigma (\delta_\Sigma \omega).$$
La même formule de Weitzenböck nous donne
\begin{eqnarray*} \Delta_\Sigma \circ \delta_\Sigma & = & \delta_\Sigma \circ
  \Delta_\Sigma \\
& = & \delta_\Sigma \circ ((\na \n)_\Sigma - (n-3)Id) \\
& = & \delta_\Sigma \circ (\na \n)_\Sigma - (n-3) \delta_\Sigma,
\end{eqnarray*}
et on a aussi $e_\theta.e_\theta.(\delta_\Sigma \omega) = \delta_\Sigma
(\n_{e_\theta} \n_{e_\theta} \omega)$, soit finalement
\begin{eqnarray*}\Delta (\delta_\Sigma \omega) & = & - \delta_\Sigma
  (\n_{e_\theta} \n_{e_\theta} \omega) +
  \frac{1}{\cosh(a)^2}\left(\delta_\Sigma((\na \n)_\Sigma \omega) -
  (n-3)\delta_\Sigma \omega\right) \\
& = & \delta_\Sigma(\na \n \omega) - \frac{n-3}{\cosh(a)^2} (\delta_\Sigma
  \omega).
\end{eqnarray*}
Comme précédemment, on constate que si $\phi_j$ est une des
sections définies plus haut, c'est-à-dire telle que $\na \n \phi_j
= \lambda_j \phi_j$ et $\n_{\frac{\d}{\d \theta}} \phi_j = i p_j
\beta \phi_j$, alors
 $$\Delta_\Sigma (\delta_\Sigma \phi_j)
=\left(\lambda_j - \frac{n-3}{\cosh(a)^2}\right) \delta_\Sigma
\phi_j$$ et que $$\frac{\d}{\d \theta} (\delta_\Sigma \phi_j) = i
p_j \beta \, \delta_\Sigma \phi_j.$$

\bigskip

Pour clarifier tout ceci, on introduit les notations suivantes. Pour $\lambda
\in \R$, $p \in \mathbb{Z}$, on pose :
\begin{eqnarray*}
E_\lambda & = &\{f \in L^2(\Sigma_a) | \Delta f = \lambda f\}, \\
F_\lambda & = &\{\omega \in L^2(N) | \na \n \omega = \lambda \omega\}, \\
E_{\lambda,p} & = &\{f \in E_\lambda | \frac{\d}{\d \theta}f = ip\beta\,f\},\\
{\rm et\ } F_{\lambda,p} & = &\{\omega \in F_\lambda | \n_{\frac{\d}{\d
    \theta}}\omega = ip\beta\,\omega\}.
\end{eqnarray*}
Chacun de ces sous-espaces vectoriels est de dimension finie,
composé de fonctions ou sections $C^\infty$. Les $E_{\lambda,p}$
sont deux à deux orthogonaux, ainsi que les $F_{\lambda,p}$, les
$E_\lambda$ et les $F_\lambda$. On a $E_\lambda =
\bigoplus_{p\in\mathbb{Z}} E_{\lambda,p}$ et $F_\lambda =
\bigoplus_{p\in\mathbb{Z}} F_{\lambda,p}$, et à chaque fois la
somme est en fait finie car $E_\lambda$ et $F_\lambda$ sont de
dimension finie. On note $S$, resp. $S'$, le spectre du laplacien
sur les fonctions, resp. sur les 1-formes, i.e. l'ensemble des
valeurs de $\lambda$ pour lesquelles $E_\lambda$, resp.
$F_\lambda$, est non réduit à $0$. Ce sont des ensembles discrets,
minorées, avec $+\infty$ comme seul point d'accumulation. Alors
$\bigoplus_{\lambda \in S} E_\lambda ( = \bigoplus_{\lambda \in S}
\bigoplus_{p\in\mathbb{Z}} E_{\lambda,p})$ est dense dans
$L^2(\Sigma_a)$, idem pour $\bigoplus_{\lambda \in S'} F_\lambda$
dans $L^2(N)$.

D'après ce qui précède, $d_\Sigma$ envoie $E_{\lambda,p}$ dans
$F_{\lambda+(n-3)/\cosh(a)^2,p}$, et $\delta_\Sigma$ envoie
$F_{\lambda+(n-3)/\cosh(a)^2,p}$ dans $E_{\lambda,p}$. De plus, si
$\psi,\psi'$ appartiennent à $E_{\lambda,p}$, alors, comme l'adjoint de
$d_\Sigma$ est $\cosh(a)^{-2}\delta_\Sigma$, en intégrant par parties on
trouve
\begin{eqnarray*}\int_{\Sigma_a} g_a(d_\Sigma \psi,d_\Sigma \psi') & = &
  \int_{\Sigma_a} \frac{1}{\cosh(a)^2}\, \psi(\delta_\Sigma d_\Sigma \psi')\\
& = &  \int_{\Sigma_a} \frac{1}{\cosh(a)^2}\, \psi\,\Delta_\Sigma \psi'\\
& = & \int_{\Sigma_a} \frac{1}{\cosh(a)^2}\, \psi\, \cosh(a)^2 (\lambda -
\frac{p_j^2 \beta^2}{\sinh(a)^2})\, \psi'\\
& = & (\lambda - \frac{p^2 \beta^2}{\sinh(a)^2}) \int_{\Sigma_a} \psi \psi'.\\
\end{eqnarray*}
En prenant $\psi=\psi'$ on trouve que, si $E_{\lambda,p}\neq \{0\}$, alors
nécessairement
$$\lambda \geq \frac{p^2 \beta^2}{\sinh(a)^2}.$$
On en déduit aussi que, pour le produit scalaire $L^2$, $d_\Sigma$ est une
homothétie de $E_{\lambda,p}$ sur son image, de rapport $\left(\lambda -
  \frac{p^2 \beta^2}{\sinh(a)^2}\right)^{1/2}$.

De plus, du fait que l'adjoint de $d_\Sigma$ est
$\cosh(a)^{-2}\delta_\Sigma$, si un élément $\phi$ de
$F_{\lambda,p}$ est dans l'orthogonal de l'image de
$E_{\lambda-(n-3)/\cosh(a)^2,p}$  par $d_\Sigma$ (ce qui est en
particulier le cas si $E_{\lambda-(n-3)/\cosh(a)^2,p} = \{0\}$),
alors nécessairement $\delta_\Sigma \phi = 0$.

On a maintenant tout ce qu'il faut pour obtenir les bases
hilbertiennes désirées. On choisit pour chaque $E_{\lambda,p}$ une
base orthonormale, leur réunion $(\psi_j)_{j\in \mathbb{N}}$ forme
une base hilbertienne de $L^2(\Sigma_a)$. Ensuite, sur chaque
$F_{\lambda,p}$, on a déjà une famille orthonormale (finie), à
savoir l'image par $(\lambda - \frac{p^2
\beta^2}{\sinh(a)^2})^{-1/2}\,d_\Sigma$ de la base orthonormale de
$E_{\lambda-(n-3)/\cosh(a)^2,p}$ (cela évidemment dans le cas où
l'on a $\lambda > \frac{p^2 \beta^2}{\sinh(a)^2}$ et
$E_{\lambda-(n-3)/\cosh(a)^2,p} \neq \{0\}$). La réunion de ces
familles orthonormales nous donne les $(\phi_j)_{j\in J}$ de la
proposition. Enfin, on complète sur chaque $F_{\lambda,p}$ cette
famille en une base orthonormée; la réunion des éléments ainsi
rajoutés constitue les $(\varphi_j)_{j\in
  \mathbb{N}}$ de la proposition.
\end{proof}

\bigskip

Notons que les éléments de ces bases hilbertiennes vérifient aussi
$$\begin{cases} \displaystyle \Delta_\Sigma \psi_j = \cosh(a)^2
  \left(\lambda_j - \frac{p_j^2
  \beta^2}{\sinh(a)^2}\right) \psi_j \\
\displaystyle (\na \n)_\Sigma \phi_j = \cosh(a)^2 \left(\lambda_j +
  \frac{n-3}{\cosh(a)^2} - \frac{p_j^2 \beta^2}{\sinh(a)^2}\right)\phi_j\\
\displaystyle (\na \n)_\Sigma \varphi_j = \cosh(a)^2 \left(\mu_j -
  \frac{{p'_j}^2\beta^2}{\sinh(a)^2}\right) \varphi_j.
\end{cases}$$

Pour simplifier ces expressions, qui sont celles qui vont nous servir, on
pose, pour tout indice $j\in \mathbb{N}$,
$$\lambda'_j = \cosh(a)^2 \left(\lambda_j - \frac{p_j^2
    \beta^2}{\sinh(a)^2}\right),$$
ainsi que
$$\mu'_j = \cosh(a)^2 \left(\mu_j -
  \frac{{p'_j}^2\beta^2}{\sinh(a)^2}\right).$$
On a alors
$$\begin{cases} \Delta_\Sigma \psi_j = \lambda'_j \psi_j \\
(\na \n)_\Sigma \phi_j = (\lambda'_j + n-3) \phi_j\\
(\na \n)_\Sigma \varphi_j = \mu'_j \varphi_j.
\end{cases}$$
On peut aussi exprimer plus simplement les relations suivantes :
$$\begin{cases}\displaystyle d_\Sigma \psi_j =
  \frac{(\lambda'_j)^{1/2}}{\cosh(a)} \phi_j \\
\delta_\Sigma \phi_j = \cosh(a) (\lambda'_j)^{1/2} \psi_j.
\end{cases}$$

\subsubsection{Expression du laplacien dans cette décomposition}

Maintenant, on va utiliser les résultats précédents pour procéder à la
décomposition de $u = f e^r + g e^\theta + \omega$ sur tout $U_a$.

Pour passer de $\Sigma_a$ à $\Sigma_r$, on utilise le transport
parallèle et le flot le long des géodésiques, intégrales du champ de vecteur
$e_r$. Cela
revient à étendre à tout $U_a$ les fonctions $\psi_j$ et les formes $\phi_j$,
$\varphi_j$, en demandant seulement que $e_r.\psi_j =0$, et que $\n_{e_r}\phi_j
= \n_{e_r}\varphi_j = 0$. On note encore $\psi_j$, $\phi_j$ et $\varphi_j$ ces
extensions.

Pour un $r$ fixé, on note $f_r$ la restriction de $f$ à
$\Sigma_r$. On peut de même étendre $f_r$ en une fonction
$\tilde{f}_r$ définie sur tout $U_a$ en utilisant le flot du champ
de vecteur $e_r$, c'est-à-dire en demandant seulement que
$e_r.\tilde{f}_r$ soit identiquement nul (et évidemment que
$\tilde{f}_r=f_r$ sur $\Sigma_r$). En particulier, on peut
regarder la restriction à $\Sigma_a$ de $\tilde{f}_r$, notée
$\tilde{f}_r|_{\Sigma_a}$ . On peut maintenant utiliser les
résultats de la proposition \ref{bashib} pour décomposer
$\tilde{f}_r|_{\Sigma_a}$ sous la forme d'une série :
$\tilde{f}_r|_{\Sigma_a} = \sum f_r^j \psi_j$. Finalement, en
réutilisant le flot pour se ramener à $\Sigma_r$, on obtient la
décomposition suivante, valable sur $\Sigma_r$ : $f_r = \sum f_r^j
\psi_j$. En faisant cette manipulation pour tout $r$, et en posant
$f_j(r) = f_r^j$, on obtient $$f = \sum_{j\in \mathbb{N}} f_j(r)
\psi_j.$$

On effectue évidemment une décomposition similaire pour la fonction $g$. Pour
la section $\omega$, le
même procédé fonctionne, en remplaçant le flot par le transport parallèle, et
on obtient une décomposition $$\omega = \sum_{j\in J} \omega_j(r) \phi_j +
\sum_{j\in \mathbb{N}} \varpi_j(r) \varphi_j.$$

On peut vérifier facilement que si $u$ est $C^\infty$ alors les coefficients
$f_j$, $g_j$, $\omega_j$ et $\varpi_j$ le sont aussi (en effet,
$f_j(r) = \int_{\Sigma_a} \overline{\psi_j}\tilde{f}_r$ et on peut dériver sous
l'intégrale; il en est de même pour les autres coefficients).

On a finalement obtenu l'expression suivante pour $u$ :
\begin{eqnarray*}
u & = & \sum_{j\in \mathbb{N}} f_j(r) \psi_j\, e^r + \sum_{j\in \mathbb{N}}
g_j(r) \psi_j\, e^\theta \\
& + & \sum_{j\in J} \omega_j(r) \phi_j + \sum_{j\in \mathbb{N}}
\varpi_j(r) \varphi_j.
\end{eqnarray*}
Il est plus judicieux de regrouper les termes de cette
décomposition de la façon suivante, faisant apparaître des ``blocs
élémentaires'' de même fréquence :
\begin{eqnarray*}
u & = & \sum_{j\in J} \left( f_j(r) \psi_j\, e^r + g_j(r) \psi_j\, e^\theta
+ \omega_j(r) \phi_j \right) \\
& + & \sum_{j \in \mathbb{N}\setminus J} \left( f_j(r) \psi_j\, e^r + g_j(r)
  \psi_j\, e^\theta \right) \\
& + & \sum_{j\in \mathbb{N}} \varpi_j(r) \varphi_j.
\end{eqnarray*}

\bigskip

En partant de cette expression pour $u$, on va effectuer cette
décomposition pour $\na \n u$. On note toujours $\beta$ pour
$\frac{2\pi}{\alpha}$. Il faut d'abord voir comment se comporte
les fonctions et formes étendues. En procédant à de simples
changement d'échelle, on arrive à :
 $$\begin{cases}\displaystyle \frac{\d}{\d r}\psi_j = 0\\
 \displaystyle \frac{\d}{\d \theta}\psi_j = ip_j\beta\,\psi_j\\
 \displaystyle \Delta_\Sigma \psi_j = \lambda'_j \psi_j\\
 \displaystyle d_\Sigma \psi_j = \frac{(\lambda'_j)^{1/2}}{\cosh(r)}
 \phi_j, \end{cases}$$
$$\begin{cases}
 \displaystyle \n_{\frac{\d}{\d r}} \phi_j = 0\\
 \displaystyle \n_{\frac{\d}{\d \theta}} \phi_j = ip_j\beta\,\phi_j \\
 \displaystyle (\na \n)_\Sigma \phi_j =  (\lambda'_j + n-3) \phi_j \\
 \displaystyle \delta_\Sigma \phi_j = \cosh(r) (\lambda'_j)^{1/2} \psi_j,
\end{cases}$$
  et
$$\begin{cases} \displaystyle \n_{\frac{\d}{\d r}} \varphi_j = 0\\
 \displaystyle \n_{\frac{\d}{\d \theta}} \varphi_j = ip'_j\beta\,\varphi_j\\
 \displaystyle (\na \n)_\Sigma \varphi_j = \mu'_j \varphi_j\\
\displaystyle  \delta_\Sigma \varphi_j = 0\end{cases}.$$

\bigskip

On obtient alors, pour la composante de $\na \n u$ en $\psi_j
e^r$, si $j \in J$ : $$ -f_j'' - \left(\frac{1}{\tanh(r)} + (n-2)
\tanh(r)\right) f_j' + \left(\frac{1}{\tanh(r)^2} + (n-2)
\tanh(r)^2 +
  \frac{p_j^2\beta^2}{\sinh(r)^2} +
\frac{\lambda'_j}{\cosh(r)^2}\right)f_j \coupeq + \frac{2ip_j\beta}{\sinh(r)
  \tanh(r)}g_j +
\frac{2\tanh(r)(\lambda'_j)^{1/2} }{\cosh(r)} \omega_j,$$

si $j \notin J$ :
$$ -f_j'' - \left(\frac{1}{\tanh(r)} + (n-2) \tanh(r)\right) f_j' +
\left(\frac{1}{\tanh(r)^2} + (n-2) \tanh(r)^2 +
  \frac{p_j^2\beta^2}{\sinh(r)^2} +
\frac{\lambda'_j}{\cosh(r)^2}\right)f_j \coupeq + \frac{2ip_j\beta}{\sinh(r)
  \tanh(r)}g_j,$$

pour la composante en $\psi_j e^\theta$ :
$$ -g_j'' - \left(\frac{1}{\tanh(r)} + (n-2) \tanh(r)\right) g_j' +
\left(\frac{1}{\tanh(r)^2}
+ \frac{p_j^2\beta^2}{\sinh(r)^2} + \frac{\lambda'_j}{\cosh(r)^2}\right)g_j
 - \frac{2ip_j\beta}{\sinh(r)\tanh(r)}f_j,$$

pour la composante en $\phi_j$ :
$$ -\omega_j'' - \left(\frac{1}{\tanh(r)} + (n-2) \tanh(r)\right)
\omega_j' + \left(\tanh(r)^2 + \frac{p_j^2\beta^2}{\sinh^2(r)} +
\frac{\lambda'_j + n-3}{\cosh(r)^2}\right)\omega_j
- \frac{2\tanh(r) (\lambda'_j)^{1/2}}{\cosh(r)} f_j,$$

et pour la composante en $\varphi_j$ :
$$ - \varpi_j'' - \left(\frac{1}{\tanh(r)} +
(n-2) \tanh(r)\right) \varpi_j' + \left(\tanh(r)^2 +
\frac{{p'_j}^2 \beta^2}{\sinh(r)^2} +
\frac{\mu'_j}{\cosh(r)^2}\right) \varpi_j.$$

\bigskip

\subsection{Comportement des solutions au voisinage de la
  singularité}\label{vois3}

On va maintenant chercher à résoudre l'équation $Lu = 0$ au voisinage de
$\Sigma$. La décomposition ci-dessus permet de passer d'une
équation aux dérivées partielles à une infinité d'équations différentielles
ordinaires. Résoudre l'équation $L u = 0$ revient donc à résoudre une équation
différentielle linéaire pour chaque coefficient de la décomposition.

L'équation qu'on obtient ici présente une singularité ``régulière" en
$r=0$. On sait (cf \cite{Wasow}) que les solutions d'une telle équation sont
des combinaisons linéaires de fonctions de la forme
$r^kf(r)$ avec $f$ une fonction analytique, où les exposants $k$ s'obtiennent
comme racines de l'équation indicielle (en cas de racines multiples ou
séparées par des entiers, il faut éventuellement
rajouter des termes en $\ln r$ dans l'expression des solutions).

On pose donc, pour un entier $j$ donné,
$$\left\{ \begin{array}{rcl}
f_j(r) & = & r^k(f_0 + f_1 r + f_2 r^2 + \cdots), \\
g_j(r) & = & r^k(g_0 + g_1 r + g_2 r^2 + \cdots), \\
\omega_j(r) & = & r^k(\omega_0 + \omega_1 r + \omega_2 r^2 + \cdots), \\
\varpi_j(r) & = & r^k(\varpi_0 +
\varpi_1 r + \varpi_2 r^2 + \cdots).
\end{array} \right.$$

On obtient alors les systèmes d'équations indicielles suivants (on omet de
noter les indices) : si $j \in J$,
$$\left\{ \begin{array}{rcccl}
(-k^2 + 1 + p^2\beta^2)f_0 & + & 2ip\beta g_0 & = & 0\\
-2ip\beta f_0 & + & (-k^2 +1 + p^2\beta^2) g_0 & = & 0\\
& &(-k^2 + p^2\beta^2)\omega_0 & = & 0, \end{array} \right.$$
si $j \notin J$,
$$\left\{ \begin{array}{rcccl}
(-k^2 + 1 + p^2\beta^2)f_0 & + & 2ip\beta g_0 & = & 0\\
-2ip\beta f_0 & + & (-k^2 +1 + p^2\beta^2) g_0 & = & 0,
\end{array} \right.$$
et enfin
$$(-k^2 + {p'}^2\beta^2)\varpi_0 = 0.$$

Commençons par étudier le premier système, le plus compliqué. Les valeurs de
l'exposant $k$ pour lesquelles il admet des solutions
non triviales (racines indicielles) sont $\pm p\beta \pm 1$ et $\pm
p\beta$. Plus précisément,
pour $k=\pm(p\beta +1)$, les coefficients dominants
$(f_0,g_0,\omega_0)$ sont engendrés par $(1,-i,0)$, pour
$k= \pm (p\beta -1)$, par $(1,i,0)$, et pour $k= \pm p\beta$, par
$(0,0,1)$. On remarque que l'on a toujours des racines
séparées par des entiers, ce qui rajoute des termes logarithmiques, mais nous
n'aurons pas à en tenir compte car seul l'exposant dominant va nous
intéresser.

Le cas des racines doubles est un peu plus compliqué. Elles apparaissent si $p
= 0$, $p\beta = \pm 1$, ou $p\beta = \pm \frac{1}{2}$. En fait si $p\beta =
\frac{1}{2}$, les solutions correspondant à $k= p\beta$ et à $k = 1 - p\beta$
sont linéairement indépendantes, on n'a donc pas besoin de termes
logarithmiques; même chose pour $ p\beta = - \frac{1}{2}$.

Pour $p\beta = 1$, les solutions pour $k = p\beta -1$ et $k = - (p\beta -1)$
sont les mêmes. On a donc besoin d'un terme logarithmique. Même chose si
$p\beta = -1$.

Enfin, pour $p = 0$, il y a trois dégénérescence. Cependant pour $k = 1$ ou $k
= -1$, on n'a pas de perte de dimension et donc pas besoin de termes
logarithmiques. Par contre, pour $k = 0$, le terme en logarithme est
nécessaire.

On remarque que les deux premiers cas de racines doubles ne se rencontrent que
pour des valeurs
particulières de l'angle conique. Par contre le dernier cas se rencontre quel
que soit l'angle. C'est l'existence de ces solutions logarithmiques, qui sont
dans $L^2$ mais dont la dérivée covariante ne l'est pas, qui fait que
l'opérateur $L$ n'est jamais essentiellement auto-adjoint dans notre cadre.

\bigskip

Les deux systèmes restant sont plus simples à étudier et ne présentent rien de
nouveau par rapport à ce qui précède. La proposition suivante regroupe tous ces
résultats :

\medskip

\begin{Prop}\label{dvpt} Soit $u$ une solution de l'équation $Lu=0$ sur un
  voisinage d'une composante connexe de $\Sigma$, d'angle conique $\alpha$.
Alors chacun des termes apparaissant dans la décomposition
\begin{eqnarray*}
u & = & \sum_{j\in J} \left( f_j(r) \psi_j\, e^r + g_j(r) \psi_j\, e^\theta
+ \omega_j(r) \phi_j \right) \\
& + & \sum_{j \in \mathbb{N}\setminus J} \left( f_j(r) \psi_j\, e^r + g_j(r)
  \psi_j\, e^\theta \right) \\
& + & \sum_{j\in \mathbb{N}} \varpi_j(r) \varphi_j.
\end{eqnarray*}
est solution de l'équation.

\bigskip

Soit $j$ un indice appartenant à $J$. L'ensemble des solutions du
type $$f_j(r) \psi_j\, e^r + g_j(r) \psi_j\, e^\theta +
\omega_j(r) \phi_j$$ forme un espace vectoriel (de dimension 6).
Si $p_j\beta \notin \{-1,0,1\}$, alors on dispose d'une base
constituée de solutions élémentaires pour lesquelles
$v(r)=(f_j(r),g_j(r),\omega_j(r))$ est de la forme $r^k(v_0+v_1
r+\cdots)$, avec $k\in \{\pm p_j \beta \pm 1, \pm p_j \beta \}$.
Pour $k=\pm(p_j\beta +1)$, on peut prendre $v_0 = (1, -i,0)$, pour
$k= \pm(p_j\beta-1)$, $(1,i,0)$, et pour $k=\pm p_j\beta$,
$(0,0,1)$. Si $p_j\beta=-1$, resp. $1$, resp. $0$, les deux
solutions élémentaires ci-dessus correspondant à $k=0$ sont
identiques, il faut donc rajouter une solution de la forme
$\ln(r)(v_0+v_1r+\cdots)$ avec $v_0=(1,-i,0)$, resp. $(1,i,0)$,
resp. $(0,0,1)$.

\bigskip

Maintenant si l'indice $j$ n'appartient pas à $J$, l'ensemble des
solutions du type $$f_j(r) \psi_j\, e^r + g_j(r) \psi_j\,
e^\theta$$ forme un espace vectoriel (de dimension 4). Si
$p_j\beta \notin \{-1,1\}$, alors on dispose d'une base constituée
de solutions élémentaires pour lesquelles $v'(r)=(f_j(r),g_j(r))$
est de la forme $r^k(v'_0+v'_1 r+\cdots)$, avec $k = \pm p_j \beta
\pm 1$. Pour $k=\pm(p_j\beta +1)$, on peut prendre $v'_0 = (1,
-i)$, et pour $k= \pm(p_j\beta-1)$, $v'_0=(1,i)$. Si
$p_j\beta=-1$, resp. $1$, les deux solutions élémentaires
ci-dessus correspondant à $k=0$ sont identiques, il faut donc
rajouter une solution de la forme $\ln(r)(v'_0+v'_1r+\cdots)$ avec
$v'_0=(1,-i)$, resp. $(1,i)$.

\bigskip

Enfin, pour tout indice $j$, l'ensemble des solutions du type
$\varpi_j(r) \varphi_j$ forme un espace vectoriel (de
dimension 2). Si $p'_j \neq 0$, alors on dispose d'une base
constituée de deux solutions élémentaires pour lesquelles
$\varpi_j(r) = r^k(1 + \varpi_1 r +\cdots)$,
avec $k = \pm p_j \beta$. Si $p'=0$ les deux solutions
élémentaires ci-dessus sont identiques,  il faut donc rajouter une
solution pour laquelle $\varpi_j(r) = \ln(r)(1 +
\varpi_1 r +\cdots)$.
\end{Prop}

\bigskip

\subsection{Résolution de l'équation}\label{2pi}

Dans cette sous-section ainsi que dans toute la suite de ce
papier, {\bf tous les angles coniques seront supposés strictement
inférieurs à $2\pi$}. En particulier, si $p$ est un entier, alors
soit $p\beta = 0$, soit $|p\beta|>1$.

On va étudier maintenant quels sont les exposants dominants
possibles pour une solution de l'équation $Lu=0$ au voisinage du
lieu singulier, en fonction des différentes conditions imposées à
$u$.

Notons tout d'abord que, vu la forme de la métrique, une forme
$u$, telle qu'au voisinage d'une composante connexe du lieu
singulier sa norme (ponctuelle) vérifie $|u| \sim r^k$, est dans
$L^2$ si et seulement si $k>-1$. Par conséquent, si $u \in L^2$
vérifie $Lu=0$ au voisinage du lieu singulier, les exposants $k$
apparaissant dans le développement de $u$ donné à la proposition
\ref{dvpt} sont tous strictement supérieurs à $-1$. Or on a vu que
$k$ est de la forme $\pm p\beta \pm 1$ ou $\pm p\beta$, où $p$ est
un entier que l'on peut supposer positif, et $\beta$ vaut $2\pi$
divisé par l'angle conique de la composante connexe du lieu
singulier. Par conséquent le fait que $u$ soit dans $L^2$ élimine
les solutions avec $k = -p\beta -1$, avec $k = -p\beta$ pour
$p\neq 0$, et avec $k = -p\beta +1$ pour $p>1$ (et aussi $p=1$ si
$\beta\geq 2$, c'est-à-dire si l'angle conique est inférieur ou
égal à $\pi$).

Le premier résultat est le lemme suivant :

\begin{Lem}\label{nnd} Soit $M$ une cône-variété hyperbolique dont tous les
  angles coniques sont strictement inférieurs à $2\pi$. Soit $u$ une 1-forme telle que
  $Lu$ soit égal à $0$ au voisinage du lieu singulier et que $u$ et $du$
  soient dans $L^2$. Alors $\n u$ et $\n du$ sont dans $L^2$.
\end{Lem}

\begin{proof} Pour $j \in J$, si $u$ est une 1-forme du type
$$f(r)\psi_j e^r + g(r)\psi_j e^\theta + \omega(r) \phi_j,$$
 alors $du$ est de la forme
 $$a(r)\psi_j e^r \wedge e^\theta + b(r) e^r
\wedge \phi_j + c(r) e^\theta \wedge \phi_j,$$
 avec (je passe les calculs)
\begin{eqnarray*}
a(r) & = & g' + \frac{1}{\tanh(r)}g -
\frac{ip_j\beta}{\sinh(r)}f,\\ b(r) & = & \omega' -\tanh(r) \omega
- \frac{(\lambda'_j)^{1/2}}{\cosh(r)} f,
\\ c(r) & = & \frac{ip_j\beta}{\sinh(r)}\omega -
\frac{(\lambda'_j)^{1/2}}{\cosh(r)} g.
\end{eqnarray*}

On suppose qu'en plus $u$ est une solution élémentaire de
l'équation $Lu=0$ au voisinage d'une composante connexe de
$\Sigma$, avec $(f(r),g(r),\omega(r))$ de la forme $r^k(v_0 + v_1
r + \cdots)$ (cf proposition \ref{dvpt}). Alors $(a(r),b(r),c(r))$
est de la forme $r^{k-1}(w_0 + w_1 r +\cdots)$, et, si on note
$v_0 = (f_0,g_0,\omega_0)$, alors
 $$w_0 = ((k+1)g_0 - ip_j\beta \,f_0,\ k
\omega_0,\ ip_j\beta\, \omega_0).$$
 On constate que si $k = \pm p_j \beta -1$, ou si $k=p_j\beta$
avec $p_j=0$, alors $w_0 = 0$, c'est-à-dire que dans ces deux cas
(et seulement dans ces deux cas-là) $u$ et $du$ ont le même
exposant dominant. Et si on est dans le cas d'un terme
logarithmique dû à une racine indicielle multiple (${p_j \beta =
-1}$, $0$, ou $1$), avec $(f(r),g(r),\omega(r))$ de la forme
$\ln(r)(v_0 + v_1 r + \cdots)$, alors l'expression de $du$
comprend toujours des termes non nuls en $r^{-1}$.

\bigskip

Maintenant pour $j \notin J$, si $u$ est une 1-forme du type
$$f(r)\psi_j e^r + g(r)\psi_j e^\theta,$$
 l'expression de $du$ est assez simple puisqu'on trouve
 $$du = \left(g' + \frac{1}{\tanh(r)}g -
\frac{ip_j\beta}{\sinh(r)} f \right) \psi_j e^r \wedge e^\theta.$$

Si $u$ est une solution élémentaire de l'équation $Lu=0$ au
voisinage (d'une composante connexe) de $\Sigma$, avec
$(f(r),g(r))$ de la forme $r^k(v'_0 + v_1 r + \cdots)$ (cf
proposition \ref{dvpt}), alors $a(r) = g' + \frac{1}{\tanh(r)}g -
\frac{ip_j\beta}{\sinh(r)} f$ est de la forme $r^{k-1}(a_0 + a_1 r
+\cdots)$, et, si on note $v'_0 = (f_0,g_0)$, alors
 $$a_0 = (k+1)g_0 - ip_j\beta \,f_0.$$
 On constate, de la même façon que dans le cas $j \in J$, que si
$k = \pm p_j \beta -1$ alors $w_0 = 0$, c'est-à-dire que dans ce
cas (et seulement dans ce cas-là) $u$ et $du$ ont le même exposant
dominant. Et si on est dans le cas d'un terme logarithmique dû à
une racine indicielle multiple (${p_j \beta = -1}$ ou $1$), avec
$(f(r),g(r))$ de la forme $\ln(r)(v'_0 + v_1 r + \cdots)$, alors
l'expression de $du$ comprend toujours des termes non nuls en
$r^{-1}$.

\bigskip

Il reste à voir ce qu'il se passe quand la solution $u$ est de la
forme $\varpi(r)\varphi_j$. On trouve, de la même
façon, que si l'exposant dominant de $u$ vaut $k$, alors
l'exposant dominant de $du$ vaut $k-1$, sauf pour la solution non
logarithmique quand $k=0$.

\bigskip

Récapitulons tout cela. Soit $u$ une solution de l'équation $Lu=0$
au voisinage d'une composante connexe de $\Sigma$, d'exposant
dominant $k$. On note $k'$ l'exposant dominant pour $du$. Alors
$k'=k-1$, sauf pour $k$ de la forme $\pm p\beta -1$ et pour $k=0$.
Donc si $u$ et $du$ sont toutes les deux dans $L^2$, alors les seules
valeurs possibles pour $k$ sont $0$, $1$, $p\beta -1$, $p\beta$ et
$p\beta+1$ (les autres valeurs pour lesquelles $u$ était $L^2$, à
savoir $k=-p\beta + 1$ et la solution logarithmique pour $k=0$,
donnent $k'\leq -1$). En particulier, on a alors $k\geq 0$ et
$k'\geq 0$.

\bigskip

Il n'est pas difficile de montrer, à l'instar de ce que l'on a
fait pour $du$, que si une solution $u$ de l'équation $Lu=0$ au
voisinage d'une composante connexe de $\Sigma$ a pour exposant
dominant $k$, alors $\n u$ a pour exposant dominant $k-1$, sauf si
$k=0$, auquel cas $\n u$ et $u$ ont le même exposant dominant
$k=0$.

Il en est de même pour $du$ et $\n du$ : si on note encore $k'$
l'exposant dominant de $du$, alors $\n du$ a pour exposant
dominant $k'-1$, sauf si $k'=0$, auquel cas $\n du$ et $du$ ont le
même exposant dominant $k'=0$.

En conclusion : si $u$ est une solution de l'équation $Lu=0$ au
voisinage de $\Sigma$, telle que $u$ et $du$ soient dans $L^2$, on
a vu que les exposants dominants $k$ et $k'$ de $u$ et de $du$
sont tous les deux supérieurs ou égaux à $0$. Par conséquent les
exposants dominants de $\n u$ et de $\n du$ sont tous les deux
strictements supérieurs à $-1$, et donc $\n u$ et $\n du$ sont
tous les deux dans $L^2$.
\end{proof}

\bigskip

L'intérêt de ce lemme réside principalement dans la démonstration des deux
résultats suivants, qui nous font passer de l'étude des solutions de
l'équation $Lu=0$ au voisinage du lieu singulier à celle des solutions de
l'équation $Lu=f$ sur $M$ entière.

\bigskip

\begin{Thm}\label{DD'} Si $M$ est une cône-variété hyperbolique dont tous les
  angles coniques sont strictement inférieurs à $2\pi$, alors $D=D'$.
\end{Thm}

\begin{proof} Supposons que les deux domaines $D$ et $D'$ soient différents;
  par exemple, $D \nsubseteq D'$. Alors il existe $\alpha \in D\backslash
  D'$. Comme $L|_{D'}$ est bijectif, il existe aussi $\alpha' \in D'$ tel que
  $L\alpha' = L\alpha$. Donc $\alpha - \alpha' \in \ker L$, et on connait le
  comportement de $ \alpha - \alpha'$ au voisinage du lieu singulier.

Par définition de $D$ et $D'$ (cf \ref{LapDeb}), on sait que $\alpha$,
$\alpha'$, $\n \alpha$ et $d\alpha'$ sont dans $L^2$, et donc aussi
$d\alpha$. Par conséquent $\alpha - \alpha'$ et $d(\alpha - \alpha')$ sont
dans $L^2$. D'après le lemme  précédent ceci implique que $\n (\alpha -
\alpha')$ est dans $L^2$.

On peut alors appliquer le théorème \ref{ipp} pour
procéder à une intégration par parties :
\begin{eqnarray*}
0 & = & \<L(\alpha-\beta),\alpha - \beta\> \\
& = & \<\na \n (\alpha-\beta) + (n-1) (\alpha-\beta),\alpha - \beta\> \\
& = & ||\n (\alpha-\beta)||^2 + (n-1) ||\alpha - \beta||^2
\end{eqnarray*}
et on trouve finalement $\alpha - \beta = 0$, ce qui contredit l'hypothèse
$\alpha \in D\backslash D'$.
\end{proof}

\medskip

{\em Remarque :} Bien que nous ne le montrions pas ici, il est
intéressant de noter que dès qu'un angle conique est plus grand
que $2\pi$, les deux domaines ci-dessus ne coïncident plus. Il
devient donc beaucoup plus difficile de trouver un ``bon'' domaine
pour résoudre l'équation de normalisation, cf \cite{HK2} et
\cite{BB} pour des résultats dans cette direction.

\bigskip

Nous allons maintenant montrer un résultat complémentaire pour les solutions
de l'équation de normalisation.

\bigskip

\begin{Thm}\label{nablad} Soit $M$ une cône-variété hyperbolique dont tous les
angles coniques sont strictement inférieurs à $2\pi$. Soit $\phi$
une section de $L^2(T^*M)$. Alors il existe une unique section
$\alpha$ de $L^2(T^*M)$, solution de l'équation $L\alpha =\phi$,
telle que $\alpha$, $\n \alpha$, $d\delta \alpha$, et $\n d\alpha$
(au sens des distributions) soient dans $L^2$.
\end{Thm}

\begin{proof} On sait depuis la section \ref{LapDeb} que l'on peut
résoudre de façon unique l'équation $L\alpha = \phi$ avec $\alpha
\in D$. Maintenant, le théorème \ref{DD'} ci-dessus nous assure
que l'on a aussi $\alpha \in D'$; finalement $\alpha$, $\n \alpha$
(et donc aussi $d \alpha$ et $\delta \alpha$), $d\delta \alpha$ et
$\delta d \alpha$ (et donc aussi $\na \n \alpha$) sont dans $L^2$.
Le seul point qui reste à montrer est que $\n d\alpha$ est aussi
$L^2$.

\bigskip

 Les formes $C^\infty$ à support compact étant dense dans $L^2$,
  on peut trouver une suite $(\phi_n)$ de 1-formes $C^\infty$ à support
  compact telle que $\phi_n \to \phi$ dans $L^2$ quand $n \to
  \infty$. Soit $(\alpha_n)$ la suite d'éléments de $D$ telle que pour tout
  entier $n$, $L \alpha_n =\phi_n$. On applique alors le théorème \ref{1+*}
  (avec, à un facteur près,
  $A=\n$ et $A^* = \na$) : les transformations $(\na \n + (n-1)Id)^{-1}$ et
  $\n(\na \n + (n-1)Id)^{-1}$ sont continues, donc
$$\lim_{n\to \infty} \alpha_n = \lim_{n\to \infty} (\na \n + (n-1)Id)^{-1}
  (\phi_n) = (\na \n + (n-1)Id)^{-1} (\phi) = \alpha$$
et
$$\lim_{n\to \infty} \n \alpha_n = \lim_{n\to \infty} \n ((\na \n +
  (n-1)Id)^{-1} (\phi_n)) = \n ((\na \n + (n-1)Id)^{-1} (\phi)) =
  \n \alpha,$$
les limites étant au sens $L^2$. Comme $d\alpha_n$ est la partie
antisymétrique de $\n \alpha_n$, la suite $(d\alpha_n)$ est aussi
convergente, avec $\lim_{n\to \infty} d \alpha_n = d\alpha$.
Maintenant, comme $\phi_n$ est à
  support compact, $L\alpha_n$ est identiquement nul au voisinage du lieu
  singulier, et $\alpha_n$ rentre donc dans le cadre de la proposition
  \ref{dvpt}. Comme $\alpha_n$ appartient à $D(=D')$, $\alpha_n$ ainsi que
  $d\alpha_n$ sont dans $L^2$, et on a vu au lemme \ref{nnd} qu'alors $\n d
  \alpha_n \in L^2$. On va maintenant montrer que $(\n d
  \alpha_n)$, suite de sections du fibré $T^*M\otimes \Lambda^2M$, est {\em
  bornée} dans $L^2$.

Pour cela, on considère $\xi$, section $C^{\infty}$ à support compact
de $T^*M\otimes \Lambda^2M$ (``section test''), et on s'intéresse au produit
scalaire $\<\n d \alpha_n, \xi \>$. Le but est d'arriver à monter que $$|\<\n d
\alpha_n, \xi \>| \leq M ||\xi||,$$ où $M$ ne dépend pas de $n$.

\bigskip

La restriction de la dérivée covariante à $\Omega^2M$ nous donne
un opérateur (non borné) $\n : L^2(\Lambda^2M) \to
L^2(T^*M\otimes\Lambda^2M)$; son adjoint est la restriction de
$\na$ à $T^*M\otimes\Lambda^2M$, et les résultats de la section
\ref{secIPP} s'appliquent. Maintenant, en utilisant la définition
de l'adjoint d'un opérateur, on a l'égalité $\ker \na = (\im
\n)^\perp$ et donc on a aussi la décomposition orthogonale
suivante :
 $$L^2(T^*M\otimes \Lambda^2M) = \ker \na \oplus
\overline{\im \n}.$$
 On voudrait pouvoir écrire $\xi = k + \n \zeta$ dans cette
décomposition, mais il faut d'abord montrer que l'image de $\n$
est fermée. Pour cela on utilise la formule de Weitzenböck
suivante, valable pour une métrique hyperbolique, qui est un
analogue de la formule pour les $1$-formes que l'on a déjà
utilisée à plusieurs reprises (cf \cite{Besse} \S 1.I) :
$$\forall \omega \in \Omega^2M,\
\na \n \omega = \Delta \omega + 2 (n-2) \omega.$$

En particulier, dès que $\omega$ est $C^\infty$ à support compact,
en intégrant par parties contre $\omega$ on obtient
 $$||\n \omega||^2 = ||d \omega||^2 + ||\delta \omega||^2 + 2 (n-2)
||\omega||^2,$$
 ce qui implique
 $$||\omega|| \leq c ||\n \omega||,$$
 avec $c=(2(n-2))^{-1/2}$. Maintenant, cette inégalité
est aussi vraie pour tout $\omega \in L^2(\Lambda^2M)$ tel que $\n
\omega \in L^2$; il suffit de prendre une suite $\omega_n \in
C^\infty_0$ telle que $\omega_n \to \omega$ et $\n \omega_n \to \n
\omega$ au sens $L^2$ (cf corollaire \ref{approx}). Cette
inégalité implique immédiatement que l'image de $\n$ est fermé, et
donc que $$L^2(T^*M\otimes \Lambda^2M) = \ker \na \oplus \im \n.$$

Donc on peut bien écrire $\xi = k + \n \zeta$, où $k \in \ker
\na$, et $||\zeta|| \leq c ||\n \zeta|| \leq c ||\xi||$.
Retournons au produit scalaire :
\begin{eqnarray*}
\<\n d \alpha_n, \xi \> & = & \<\n d \alpha_n, \n \zeta + k\> \\
& = & \<\n d \alpha_n, \n \zeta\>.
\end{eqnarray*}
Pour pouvoir faire une intégration par parties, il faut vérifier
que tous les termes impliqués sont $L^2$. On sait déjà que
$\zeta$, $\n \zeta$, $\n d \alpha_n$ le sont. Or d'après la
formule de Weitzenböck ci-dessus,
 $$\na \n d \alpha_n = \Delta d
\alpha_n + 2(n-2) d \alpha_n = d \Delta \alpha_n + 2(n-2) d
\alpha_n,$$
 car $d$ et $\Delta = d\delta +\delta d$ commutent.
D'autre part
 $$\Delta \alpha_n = L\alpha_n - 2(n-1)\alpha_n =
\phi_n - 2 (n-1) \alpha_n.$$
 Finalement,
$$\na \n d \alpha_n = d \phi_n - 2 d\alpha_n.$$

Comme $\phi_n$ est à support compact et que $\alpha_n \in D$, les
formes $d\phi_n$ et $d\alpha_n$ sont $L^2$, donc $\na \n d
\alpha_n$ est $L^2$, donc on peut donc intégrer par parties
(théorème \ref{ipp}) :
\begin{eqnarray*}
\<\n d \alpha_n, \xi \> & = & \<\n d \alpha_n, \n \zeta\>\\ & = &
\<\na \n d \alpha_n, \zeta\>\\ & = & \<d \phi_n - 2 d\alpha_n,
\zeta\>.
\end{eqnarray*}

Comme $\n \zeta$ est $L^2$, $\delta \zeta = - \tr_g \n \zeta$ est
aussi $L^2$, on a même $||\delta \zeta|| \leq \sqrt{n} ||\n
\zeta||$. D'autre part $\phi_n$, $d\phi_n$ et $\zeta$ sont $L^2$,
on peut encore intégrer par parties :
\begin{eqnarray*}
\<\n d \alpha_n, \xi \> & = & \<d \phi_n - 2 d\alpha_n, \zeta\>\\
& = & \<\phi_n,\delta \zeta\>- 2\<d\alpha_n,\zeta\>.
\end{eqnarray*}

Pour finir on majore avec Cauchy-Schwarz:
\begin{eqnarray*}
|\<\n d \alpha_n, \xi \>| & \leq & ||\phi_n||\,||\delta \zeta|| +
2 ||d\alpha_n||\,||\zeta|| \\ & \leq & (\sqrt{n} ||\phi_n|| + 2 c
||d\alpha_n||)\,||\n \zeta||\\ & \leq & M ||\xi||
\end{eqnarray*}
car les suites $(\phi_n)$ et $(d\alpha_n)$ sont convergentes, donc bornées,
dans $L^2$. Cette majoration, valable pour toute section test $\xi$, implique
directement que la suite $(\n d\alpha_n)$ est bornée dans $L^2$.

\bigskip

Par conséquent, on peut extraire une sous-suite, encore notée $(\n
d\alpha_n)$, qui converge faiblement vers une limite $l\in L^2$ : c'est-à-dire
que quel que soit $\xi \in L^2(T^*M\otimes \Lambda^2M)$,
$$\lim_{n \to \infty}\<\n d\alpha_n, \xi\> = \<l,\xi\>.$$ Mais alors, si $\xi$
est $C^\infty$ à
support compact,
$$\<\n d\alpha_n, \xi\> = \<d\alpha_n, \na \xi\>,$$ et
$$\lim_{n \to \infty} \<d\alpha_n, \na \xi\> = \<d\alpha,\na \xi\>$$ car
$(d\alpha_n)$ converge dans $L^2$ vers $d\alpha$. Par conséquent, on a
$$\<d\alpha,\na \xi\> = \<l,\xi\>$$ pour tout $\xi \in C^\infty_0$, ce qui
signifie exactement que $$l=\n_{max} d\alpha = \n d\alpha,$$
 et par suite $\n d\alpha$ appartient à $L^2$.
\end{proof}

\medskip

Notons que si en plus $\phi$ est $C^\infty$, alors par régularité
elliptique la solution $\alpha$ ci-dessus est aussi de classe
$C^\infty$.

\bigskip

\section{Rigidité infinitésimale des c\^one-variétés}\label{rigidite}

Nous avons maintenant en main tous les outils pour montrer le théorème suivant
:

\begin{Thm}\label{thmppal} Soit $M$ une c\^one-variété hyperbolique dont tous
  les angles coniques sont strictement inférieurs à $2\pi$.
Soit $h_0$ une déformation Einstein infinitésimale (i.e. vérifiant
  l'équation $E'_g(h_0)=0$) telle que $h_0$ et $\n h_0$ soient dans
  $L^2$. Alors la déformation $h_0$ est triviale, i.e. il existe une forme
  $\alpha \in \Omega^1 M$ telle que $h_0 = \delta^* \alpha$.
\end{Thm}

Dans toute cette section nous supposerons donc que les angles coniques sont
toujours inférieurs à $2\pi$.

\bigskip

\begin{proof}
La première étape de la démonstration consiste à normaliser $h_0$,
c'est-à-dire à chercher $\alpha$ tel que $h = h_0 - \delta^*
\alpha$ vérifie la condition de jauge $\beta(h)=0$, ce qui revient
à résoudre l'équation $\beta \circ \delta^* \alpha = \beta h_0$.
Comme $\n h_0$ est dans $L^2$, $\beta h_0$ l'est aussi, et d'après
le théorème \ref{nablad} cette équation admet une unique solution
$\alpha$ telle que $\alpha$, $\n \alpha$, $d \delta \alpha$ et $\n
d\alpha$ soient dans $L^2$. On pose $h = h_0 - \delta^*\alpha$.
Notons que l'on a perdu des informations en normalisant : en
effet, rien ne garantit que la déformation normalisée $h$ vérifie
encore $\n h =0$, puisqu'on ne connaît rien pour l'instant sur $\n
\delta^* \alpha$.

La déformation $h$ vérifie alors :
$$\begin{cases} \na \n h - 2 \mathring{R}h = 0 \\ \delta h + d \tr h = 0
\end{cases}$$

En prenant la trace par rapport à $g$ de la première équation, on
obtient $$\Delta(\tr h) + 2(n-1) \tr h =0,$$ ce qui incite à
intégrer par parties, mais pour le faire il faut d'abord vérifier
que les termes impliqués sont $L^2$, avant de pouvoir appliquer le
théorème \ref{Cheeg}. Comme $h_0$ et $\delta^* \alpha$ sont $L^2$,
$h$ est bien $L^2$, donc $\tr h$ aussi, et donc $\Delta \tr h$
aussi. Maintenant, $$d \tr h = d \tr h_0 + d \tr \delta^* \alpha =
d \tr h_0 - d\delta\alpha,$$
 donc $d \tr h$ est $L^2$ ($d \tr h_0$ est $L^2$ car $\n h_o$ l'est). Par
suite, on trouve en intégrant contre $\tr h$ :
\begin{eqnarray*} 0 & = & \< \tr h, \Delta (\tr h) + 2 (n-1)
\tr h \> \\
 & = & ||d\tr h||^2 + 2(n-1)||\tr h||^2
\end{eqnarray*}
et donc $\tr h = 0$, ce qui, avec $\beta(h)=0$, implique aussi
$\delta h = 0$. Finalement, on a
 $$\begin{cases} \na \n h - 2
\mathring{R}h = 0 \\ \delta h =0 \\ \tr h = 0
\end{cases}$$

\bigskip

La deuxième étape de la démonstration consiste à utiliser une
autre formule de Weitzenböck (cf \cite{Besse}, \S 12.69). Un
$2$-tenseur peut toujours se voir comme une $1$-forme à valeur
dans le fibré cotangent $T^*M$. Ce fibré étant muni de la
connexion de Levi-Cività $\n$, on note $d^\n$ la différentielle
extérieure associée sur les formes à valeurs dans $T^*M$.
L'opérateur adjoint est la codifférentielle notée $\delta^\n$.
Notons que si $\alpha$ est une $0$-forme à valeurs dans $T^*M$
(c'est-à-dire une $1$-forme usuelle), alors $d^\n \alpha = \n
\alpha$; de même pour une $1$-forme à valeurs dans $T^*M$,
$\delta^\n h = \na h$. On a alors la formule suivante, valable
pour tout $2$-tenseur symétrique :
 $$\na \n h = (\delta^\n d^\n + d^\n \delta^\n) h + \mathring{R}h - h\circ
 ric.$$
Pour une métrique hyperbolique, cela se simplifie en
$$\na \n h = (\delta^\n d^\n + d^\n \delta^\n) h + n h - (\tr h)g.$$
En combinant avec ce qui précède, on obtient
$$\begin{cases}  \delta^\n d^\n h + (n-2) h = 0 \\ \delta h =0 \\ \tr h = 0
\end{cases}$$

Pour conclure, ``il suffit" d'une intégration par parties contre $h$. Comme $h$
est dans
$L^2$, $\delta^\n d^\n h$ est aussi dans $L^2$; si $\n h$, ou même seulement
$\n_{e_r} h$, était $L^2$ on pourrait conclure en utilisant une méthode
analogue à celle employée dans la démonstration du théorème
\ref{ipp}. Malheureusement on ne sait rien sur le caractère $L^2$ ou non de
$\n \delta^* \alpha$. On va donc devoir contourner cette difficulté pour
montrer qu'on a bien $\<\delta^\n d^\n h, h\> = ||d^\n h||^2$.

Avant toutes choses, il faut montrer que $d^\n h$ est bien $L^2$. Comme $\n
h_0$ est $L^2$, $d^\n h_0$ est $L^2$; il ne reste qu'à regarder $d^\n
\delta^*\alpha$. Or
$$\delta^*\alpha = \n \alpha - \frac{1}{2}d\alpha = d^\n
\alpha - \frac{1}{2}d\alpha,$$
 donc
$$d^\n \delta^*\alpha = (d^\n)^2\alpha -
\frac{1}{2} d^\n d\alpha.$$
L'opérateur $(d^\n)^2$ est bien connu, ce n'est rien d'autre
que l'opposé de la courbure, i.e.
$$(d^\n)^2\alpha (x,y) = - R(x,y)\alpha = \n_x \n_y \alpha - \n_y \n_x \alpha
 - \n_{[x,y]}\alpha.$$
 C'est un opérateur borné, c'est-à-dire continue, pour les normes $L^2$; par
 conséquent
$(d^\n)^2\alpha$ est $L^2$. Il ne nous reste donc que la terme $d^\n d\alpha$;
or le théorème \ref{nablad} nous garantit que $\n d\alpha$, et donc $d^\n
d\alpha$, sont bien $L^2$.

\bigskip

Le tenseur $d^\n h$ est donc bien dans $L^2$. Malheureusement,
on n'a pas d'analogue du résultat de Cheeger (théorème \ref{Cheeg}) pour les
formes à valeurs dans un fibré, du fait que $d^\n\circ d^\n$ ne s'annule pas
nécessairement, à la différence de $d\circ d$. Cependant, en écrivant
$$h = h_0 - \delta^*\alpha = h_0 + \frac{1}{2}d\alpha - d^\n \alpha,$$
 on a
$$\<h,\delta^\n d^\n h\> = \<h_0 + \frac{1}{2}d\alpha, \delta^\n d^\n h\> -
\<d^\n \alpha, \delta^\n d^\n h\>.$$

Le théorème \ref{nablad} nous assure que $\n(h_0 + \frac{1}{2}d\alpha)$ est
dans $L^2$. Ceci
nous permet de montrer, exactement de la même façon que dans la démonstration
du théorème \ref{ipp}, qu'on a bien
$$\<h_0 + \frac{1}{2}d\alpha, \delta^\n d^\n h\> = \<d^\n(h_0 +
\frac{1}{2}d\alpha), d^\n h\>.$$

Pour le terme qui reste, comme $\alpha$ et $\n \alpha$ sont $L^2$, on peut
trouver d'après le corollaire \ref{approx} une suite $(\alpha_n)$, $C^\infty$
à support compact, telle que $\lim_{n\to\infty} \alpha_n = \alpha$ et
$\lim_{n\to\infty} \n \alpha_n = \n \alpha$. On a alors
$$\lim_{n\to \infty} \<d^\n \alpha_n, \delta^\n d^\n h\> =
\<d^\n \alpha, \delta^\n d^\n h\>.$$
 On peut faire l'intégration par parties avec
$\alpha_n$ :
$$\<d^\n \alpha_n, \delta^\n d^\n h\> = \<(d^\n)^2 \alpha_n, d^\n
h\>.$$
 Mais comme $(d^\n)^2$ est continue, on a
$$\lim_{n\to \infty} (d^\n)^2 \alpha_n = (d^\n)^2 \alpha,$$
 et donc
$$\lim_{n\to \infty} \<(d^\n)^2 \alpha_n, d^\n h\> = {\<(d^\n)^2 \alpha, d^\n
 h\>}.$$
 On en déduit que
$$\<d^\n \alpha, \delta^\n d^\n h\> = \<(d^\n)^2 \alpha, d^\n h\>,$$
et avec ce qui précède on a établi l'égalité
$$\<h,\delta^\n d^\n h\> = ||d^\n h||^2.$$

Par conséquent, comme $\delta^\n d^\n h + (n-2) h =0$, on a
\begin{eqnarray*} 0 & = & \<h,\delta^\n d^\n h + (n-2) h\> \\
& = & ||d^\n h||^2 + (n-2) ||h||^2
\end{eqnarray*}
et donc le tenseur $h$ est identiquement nul. Par suite $h_0 = \delta^*
\alpha$, la déformation est triviale.
\end{proof}

\bigskip

\begin{Cor} Soit $M$ une c\^one-variété hyperbolique dont tous
  les angles coniques sont strictement inférieurs à $2\pi$. Alors $M$ est
  infinitésimalement rigide parmi les cônes-variétés Einstein à angles
  coniques fixés.
\end{Cor}

\begin{proof} En effet, on a vu que toute déformation infinitésimale de la
  structure de cône-variété préservant les angles pouvait se mettre sous la
  forme d'un $2$-tenseur symétrique $h_0$ appartenant à $L^2$, dont la dérivée
  covariante $\n h_0$ est aussi dans $L^2$. On peut alors appliquer le théorème
  ci-dessus pour montrer que toutes les déformations Einstein de ce type sont
  triviales.
\end{proof}


\end{document}